\documentclass[12pt]{amsart}
\usepackage{amssymb,latexsym,tikz}
\usepackage[top=40mm, bottom=25mm, left=35mm, right=35mm]{geometry}
\usetikzlibrary{arrows,decorations.markings}
\tikzstyle arrowstyle=[scale=1]
\tikzstyle directed=[postaction={decorate,
decoration={markings,mark=at position .65 with {\arrow[arrowstyle]{stealth}}}}]

\usepackage{accents}
\usepackage{mathtools}
\usepackage{tikz-cd}
\usetikzlibrary{decorations.markings}
\tikzset{mid vert/.style={/utils/exec=\tikzset{every node/.append style={outer sep=0.8ex}},
postaction=decorate,decoration={markings,
mark=at position 0.5 with {\draw[-] (0,#1) -- (0,-#1);}}},
mid vert/.default=0.75ex}
\usepackage{amscd}
\usepackage{xypic}
\usepackage[utf8]{inputenc}
\usepackage{graphics}

\begin{document}

\title 
{Profunctors between posets and Alexander duality}
       
\author{Gunnar Fl{\o}ystad}
\address{Matematisk Institutt\\
         Postboks\\
         5020 Bergen}
\email{nmagf@uib.no}


\date{\today}


\theoremstyle{plain}
\newtheorem{theorem}{Theorem}[section]
\newtheorem{corollary}[theorem]{Corollary}
\newtheorem*{main}{Main Theorem}
\newtheorem{lemma}[theorem]{Lemma}
\newtheorem{proposition}[theorem]{Proposition}
\newtheorem{conjecture}[theorem]{Conjecture}

\theoremstyle{definition}
\newtheorem{definition}[theorem]{Definition}
\newtheorem{fact}[theorem]{Fact}
\newtheorem{obs}[theorem]{Observation}
\newtheorem{problem}[theorem]{Problem}

\theoremstyle{remark}
\newtheorem{notation}[theorem]{Notation}
\newtheorem{remark}[theorem]{Remark}
\newtheorem{example}[theorem]{Example}
\newtheorem{condition}[theorem]{Condition}
\newtheorem{claim}{Claim}
\newtheorem{observation}[theorem]{Observation}


\newcommand{\psp}[1]{{{\bf P}^{#1}}}
\newcommand{\psr}[1]{{\bf P}(#1)}
\newcommand{\op}{{\mathcal O}}
\newcommand{\opw}{\op_{\psr{W}}}

\newcommand{\ini}[1]{\text{in}(#1)}
\newcommand{\gin}[1]{\text{gin}(#1)}
\newcommand{\kr}{{\Bbbk}}
\newcommand{\pd}{\partial}
\newcommand{\vardel}{\partial}
\renewcommand{\tt}{{\bf t}}


\newcommand{\coh}{{{\text{{\rm coh}}}}}


\newcommand{\modv}[1]{{#1}\text{-{mod}}}
\newcommand{\modstab}[1]{{#1}-\underline{\text{mod}}}

\newcommand{\sut}{{}^{\tau}}
\newcommand{\sumit}{{}^{-\tau}}
\newcommand{\til}{\thicksim}

\newcommand{\totp}{\text{Tot}^{\prod}}
\newcommand{\dsum}{\bigoplus}
\newcommand{\dprod}{\prod}
\newcommand{\lsum}{\oplus}
\newcommand{\lprod}{\Pi}

\newcommand{\La}{{\Lambda}}

\newcommand{\sirstj}{\circledast}

\newcommand{\she}{\EuScript{S}\text{h}}
\newcommand{\cm}{\EuScript{CM}}
\newcommand{\cmd}{\EuScript{CM}^\dagger}
\newcommand{\cmri}{\EuScript{CM}^\circ}
\newcommand{\cler}{\EuScript{CL}}
\newcommand{\clerd}{\EuScript{CL}^\dagger}
\newcommand{\clerri}{\EuScript{CL}^\circ}
\newcommand{\gor}{\EuScript{G}}
\newcommand{\cF}{\mathcal{F}}
\newcommand{\cG}{\mathcal{G}}
\newcommand{\cM}{\mathcal{M}}
\newcommand{\cE}{\mathcal{E}}
\newcommand{\cI}{\mathcal{I}}
\newcommand{\cP}{\mathcal{P}}
\newcommand{\cK}{\mathcal{K}}
\newcommand{\cS}{\mathcal{S}}
\newcommand{\cC}{\mathcal{C}}
\newcommand{\cO}{\mathcal{O}}
\newcommand{\cJ}{\mathcal{J}}
\newcommand{\cU}{\mathcal{U}}
\newcommand{\cQ}{\mathcal{Q}}
\newcommand{\cX}{\mathcal{X}}
\newcommand{\cY}{\mathcal{Y}}
\newcommand{\cZ}{\mathcal{Z}}
\newcommand{\cV}{\mathcal{V}}

\newcommand{\mm}{\mathfrak{m}}

\newcommand{\dlim} {\varinjlim}
\newcommand{\ilim} {\varprojlim}

\newcommand{\CM}{\text{CM}}
\newcommand{\Mon}{\text{Mon}}


\newcommand{\Kom}{\text{Kom}}


\newcommand{\EH}{{\mathbf H}}
\newcommand{\res}{\text{res}}
\newcommand{\Hom}{\text{Hom}}
\newcommand{\inhom}{{\underline{\text{Hom}}}}
\newcommand{\Ext}{\text{Ext}}
\newcommand{\Tor}{\text{Tor}}
\newcommand{\ghom}{\mathcal{H}om}
\newcommand{\gext}{\mathcal{E}xt}
\newcommand{\id}{\text{{id}}}
\newcommand{\im}{\text{im}\,}
\newcommand{\codim} {\text{codim}\,}
\newcommand{\resol}{\text{resol}\,}
\newcommand{\rank}{\text{rank}\,}
\newcommand{\lpd}{\text{lpd}\,}
\newcommand{\coker}{\text{coker}\,}
\newcommand{\supp}{\text{supp}\,}
\newcommand{\Ad}{A_\cdot}
\newcommand{\Bd}{B_\cdot}
\newcommand{\Fd}{F_\cdot}
\newcommand{\Gd}{G_\cdot}


\newcommand{\sus}{\subseteq}
\newcommand{\sups}{\supseteq}
\newcommand{\pil}{\rightarrow}
\newcommand{\vpil}{\leftarrow}
\newcommand{\rpil}{\leftarrow}
\newcommand{\lpil}{\longrightarrow}
\newcommand{\inpil}{\hookrightarrow}
\newcommand{\pils}{\twoheadrightarrow}
\newcommand{\projpil}{\dashrightarrow}
\newcommand{\dotpil}{\dashrightarrow}
\newcommand{\adj}[2]{\overset{#1}{\underset{#2}{\rightleftarrows}}}
\newcommand{\mto}[1]{\stackrel{#1}\longrightarrow}
\newcommand{\vmto}[1]{\stackrel{#1}\longleftarrow}
\newcommand{\mtoelm}[1]{\stackrel{#1}\mapsto}
\newcommand{\bihom}[2]{\overset{#1}{\underset{#2}{\rightleftarrows}}}
\newcommand{\eqv}{\Leftrightarrow}
\newcommand{\impl}{\Rightarrow}

\newcommand{\iso}{\cong}
\newcommand{\te}{\otimes}
\newcommand{\into}[1]{\hookrightarrow{#1}}
\newcommand{\ekv}{\Leftrightarrow}
\newcommand{\equi}{\simeq}
\newcommand{\isopil}{\overset{\cong}{\lpil}}
\newcommand{\equipil}{\overset{\equi}{\lpil}}
\newcommand{\ispil}{\isopil}
\newcommand{\vvi}{\langle}
\newcommand{\hvi}{\rangle}
\newcommand{\susneq}{\subsetneq}
\newcommand{\sgn}{\text{sign}}


\newcommand{\xd}{\check{x}}
\newcommand{\ortog}{\bot}
\newcommand{\tL}{\tilde{L}}
\newcommand{\tM}{\tilde{M}}
\newcommand{\tH}{\tilde{H}}
\newcommand{\tvH}{\widetilde{H}}
\newcommand{\tvh}{\widetilde{h}}
\newcommand{\tV}{\tilde{V}}
\newcommand{\tS}{\tilde{S}}
\newcommand{\tT}{\tilde{T}}
\newcommand{\tR}{\tilde{R}}
\newcommand{\tf}{\tilde{f}}
\newcommand{\ts}{\tilde{s}}
\newcommand{\tp}{\tilde{p}}
\newcommand{\tr}{\tilde{r}}
\newcommand{\tfst}{\tilde{f}_*}
\newcommand{\empt}{\emptyset}
\newcommand{\bfa}{{\mathbf a}}
\newcommand{\bfb}{{\mathbf b}}
\newcommand{\bfd}{{\mathbf d}}
\newcommand{\bfl}{{\mathbf \ell}}
\newcommand{\bfx}{{\mathbf x}}
\newcommand{\bfm}{{\mathbf m}}
\newcommand{\bfv}{{\mathbf v}}
\newcommand{\bft}{{\mathbf t}}
\newcommand{\bbfa}{{\mathbf a}^\prime}
\newcommand{\la}{\lambda}
\newcommand{\bfen}{{\mathbf 1}}
\newcommand{\bfe}{{\mathbf 1}}
\newcommand{\ep}{\epsilon}
\newcommand{\en}{r}
\newcommand{\tu}{s}
\newcommand{\Sym}{\text{Sym}}

\newcommand{\ome}{\omega_E}

\newcommand{\bevis}{{\bf Proof. }}
\newcommand{\demofin}{\qed \vskip 3.5mm}
\newcommand{\nyp}[1]{\noindent {\bf (#1)}}
\newcommand{\demo}{{\it Proof. }}
\newcommand{\demodone}{\demofin}
\newcommand{\parg}{{\vskip 2mm \addtocounter{theorem}{1}  
                   \noindent {\bf \thetheorem .} \hskip 1.5mm }}

\newcommand{\lcm}{{\text{lcm}}}


\newcommand{\dl}{\Delta}
\newcommand{\cdel}{{C\Delta}}
\newcommand{\cdelp}{{C\Delta^{\prime}}}
\newcommand{\dlst}{\Delta^*}
\newcommand{\Sdl}{{\mathcal S}_{\dl}}
\newcommand{\lk}{\text{lk}}
\newcommand{\lkd}{\lk_\Delta}
\newcommand{\lkp}[2]{\lk_{#1} {#2}}
\newcommand{\del}{\Delta}
\newcommand{\delr}{\Delta_{-R}}
\newcommand{\dd}{{\dim \del}}
\newcommand{\Del}{\Delta}

\renewcommand{\aa}{{\bf a}}
\newcommand{\bb}{{\bf b}}
\newcommand{\cc}{{\bf c}}
\newcommand{\xx}{{\bf x}}
\newcommand{\yy}{{\bf y}}
\newcommand{\zz}{{\bf z}}
\newcommand{\mv}{{\xx^{\aa_v}}}
\newcommand{\mF}{{\xx^{\aa_F}}}

\newcommand{\Symm}{\text{Sym}}
\newcommand{\pnm}{{\bf P}^{n-1}}
\newcommand{\opnm}{{\go_{\pnm}}}
\newcommand{\ompnm}{\omega_{\pnm}}

\newcommand{\pn}{{\bf P}^n}
\newcommand{\hele}{{\mathbb Z}}
\newcommand{\nat}{{\mathbb N}}
\newcommand{\rasj}{{\mathbb Q}}
\newcommand{\bfone}{{\mathbf 1}}

\newcommand{\dt}{\bullet}
\newcommand{\disk}{\scriptscriptstyle{\bullet}}

\newcommand{\cxF}{F_\dt}
\newcommand{\pol}{f}

\newcommand{\Rn}{{\mathbb R}^n}
\newcommand{\An}{{\mathbb A}^n}
\newcommand{\frg}{\mathfrak{g}}
\newcommand{\PW}{{\mathbb P}(W)}

\newcommand{\pos}{{\mathcal Pos}}
\newcommand{\g}{{\gamma}}

\newcommand{\Vaa}{V_0}
\newcommand{\Bp}{B^\prime}
\newcommand{\Bpp}{B^{\prime \prime}}
\newcommand{\bbp}{\mathbf{b}^\prime}
\newcommand{\bbpp}{\mathbf{b}^{\prime \prime}}
\newcommand{\bp}{{b}^\prime}
\newcommand{\bpp}{{b}^{\prime \prime}}

\newcommand{\oLa}{\overline{\Lambda}}
\newcommand{\ov}[1]{\overline{#1}}
\newcommand{\ovv}[1]{\overline{\overline{#1}}}
\newcommand{\tm}{\tilde{m}}
\newcommand{\po}{\bullet}

\newcommand{\surj}[1]{\overset{#1}{\twoheadrightarrow}}
\newcommand{\Supp}{\text{Supp}}

\def\CC{{\mathbb C}}
\def\GG{{\mathbb G}}
\def\ZZ{{\mathbb Z}}
\def\NN{{\mathbb N}}
\def\RR{{\mathbb R}}
\def\OO{{\mathbb O}}
\def\QQ{{\mathbb Q}}
\def\VV{{\mathbb V}}
\def\PP{{\mathbb P}}
\def\EE{{\mathbb E}}
\def\FF{{\mathbb F}}
\def\AA{{\mathbb A}}

\newcommand{\oR}{\overline{R}}
\newcommand{\bfu}{{\mathbf u}}
\newcommand{\nn}{{\mathbf n}}
\newcommand{\oa}{\overline{a}}
\newcommand{\cop}{\text{cop}}
\renewcommand{\op}{\text{op}\,}
\renewcommand{\mm}{{\mathbf m}}
\newcommand{\ngmi}{\text{neg}}
\newcommand{\up}{\text{up}}
\newcommand{\dw}{\text{down}}
\newcommand{\diw}[1]{\widehat{#1}}
\newcommand{\di}[1]{\diw{#1}}
\newcommand{\bo}{b}
\newcommand{\ub}{u}
\newcommand{\fs}{\infty}
\newcommand{\ifst}{\infty}
\newcommand{\mon}{{mon}}
\newcommand{\cl}{\text{cl}}
\newcommand{\rego}{\text{reg\,}}
\newcommand{\ul}[1]{\underline{#1}}
\renewcommand{\ov}[1]{\overline{#1}}
\newcommand{\bipil}{\leftrightarrow}
\newcommand{\bfc}{{\mathbf c}}
\renewcommand{\mp}{m^\prime}
\newcommand{\np}{n^\prime}
\newcommand{\Mod}{\text{Mod }}
\newcommand{\Sh}{\text{Sh } }
\newcommand{\st}{\text{st}}
\newcommand{\hM}{\tilde{M}}
\newcommand{\hs}{\tilde{s}}
\newcommand{\ee}{\mathbf{e}}
\renewcommand{\dd}{\mathbf{d}}
\renewcommand{\en}{{\mathbf 1}}
\long\def\ignore#1{}
\newcommand{\lex}{{\text{lex}}}
\newcommand{\ordGL}{\succeq_{\lex}}
\newcommand{\ordG}{\succ_{\lex}}
\newcommand{\ordML}{\preceq_{\lex}}
\newcommand{\ordM}{\prec_{\lex}}
\newcommand{\tLa}{\tilde{\Lambda}}
\newcommand{\tGa}{\tilde{\Gamma}}
\newcommand{\STS}{\text{STS}}
\newcommand{\ii}{{\rotatebox[origin=c]{180}{\scalebox{0.7}{\rm{!}}}}}
\newcommand{\jj}{\mathbf{j}}
\renewcommand{\mod}{\text{ mod}\,}
\newcommand{\shmod}{\texttt{shmod}\,}
\newcommand{\hf}{\underline{f}}
\newcommand{\Glim}{\lim}
\newcommand{\Gcolim}{\colim}
\newcommand{\fm}{f^{\underline{m}}}
\newcommand{\fn}{f^{\underline{n}}}
\newcommand{\gn}{g_{\mathbf{|}n}}
\newcommand{\se}[1]{\overline{#1}}
\newcommand{\cB}{{\mathcal{B}}}
\newcommand{\fin}{{\text{fin}}}
\newcommand{\Poset}{{\text{\bf Poset}}}
\newcommand{\Set}{{\text{\bf Set}}}
\newcommand{\Homi}{\text{Hom}^{L}}
\newcommand{\Homb}{\text{Hom}_S}
\newcommand{\Homu}{\text{Hom}^u}
\newcommand{\bfnu}{{\mathbf 0}}
\renewcommand{\bfen}{{\mathbf 1}}
\newcommand{\semip}{up semi-finite }
\newcommand{\semim}{down semi-finite }
\newcommand{\difi}[1]{\di{#1}_{fin}}
\newcommand{\diff}[1]{\di{#1}^{fin}}
\newcommand{\ddfin}[1]{\doublehat{#1}\!-\!{\rm fin}}
\newcommand{\llin}{\raisebox{1pt}{\scalebox{1}[0.6]{$\mid$}}}
\newcommand{\promap}{\mathrlap{{\hskip 2.8mm}{\llin}}{\lpil}}
\newcommand{\dual}{{\widehat{}}}
\newcommand{\Bool}{{\bf Bool}}
\newcommand{\pro}{{pro}}

\newlength{\dhatheight}
\newcommand{\doublehat}[1]{%
    \settoheight{\dhatheight}{\ensuremath{\widehat{#1}}}%
    \addtolength{\dhatheight}{-0.2ex}
    \widehat{\vphantom{\rule{10pt}{\dhatheight}}
    \smash{\widehat{#1}}}}
\newcommand{\ddhat}[1]{\widehat{\!\widehat{#1}}}
\newcommand{\dfhat}[1]{\overset{-}{\widehat{#1}}}
\newcommand{\ddbhat}[1]{\widehat{\widehat{#1}}}

\newcommand{\napo}{natural }

\newcommand{\colim}{\text{colim}}
\newcommand{\upa}{\uparrow}
\newcommand{\doa}{\downarrow}
\newcommand{\boto}{{\mathbf 2}}
\newcommand{\Pro}{\text{Pro}}
\newcommand{\Proi}{\text{Pro}^L}
\newcommand{\Prob}{\text{Pro}_S}
\newcommand{\Prou}{\text{Pro}^u}
\newcommand{\DP}{{\mathbb D} P}

\begin{abstract}
  We consider profunctors $f : P \promap Q$ between posets
  and introduce their {\em graph} and {\em ascent}.
  The profunctors $\Pro(P,Q)$ form themselves a poset, and we consider
  a partition $\cI \sqcup \cF$ of this into a down-set $\cI$ and up-set $\cF$,
  called a {\it cut}.  To elements of $\cF$ we associate their graphs, 
  and to elements of $\cI$ we associate their ascents. Our basic results
  is that this, suitable refined, preserves being a cut:
  We get a cut in the Boolean lattice of subsets of the underlying set of
  $Q \times P$.
  Cuts in finite Booleans lattices correspond precisely to finite simplicial
  complexes.

  We apply this in commutative algebra where these
  give classes of Alexander dual
  square-free monomial ideals giving the full and natural generalized setting of
  isotonian ideals and
  letterplace ideals for
  posets. We study $\Pro(\NN, \NN)$.
  Such profunctors identify as
  order preserving maps
  $f : \NN \pil \NN \cup \{\infty \}$. For our applications when $P$ and
  $Q$ are infinite, we also introduce a topology on
  $\Pro(P,Q)$, in particular on profunctors
  $\Pro(\NN,\NN)$.
\end{abstract}
\maketitle

%
%
%
%
%

\section*{Introduction}

This article advocates for general posets $P$ and $Q$ the notion of 
profunctor $P \promap Q$ as more effective
than the notion of isotone (order preserving) maps $P \pil Q$
between posets, especially for applications in algebra.
When $Q$ is totally ordered, these notions are practically
the same, but when $Q$ is not, profunctors seem to have
a clear advantage for developing natural theory.

Let $\boto$ be the two element boolean poset $\{0 < 1 \}$.
A profunctor $P \promap Q$ is simply
an isotone map $P \times Q^\op \pil \boto$. If $P$ and $Q$ are sets
(discrete posets), then this is simply a relation between $P$ and $Q$.
The notion of profunctor
may generally be defined between categories or
between categories enriched in a symmetric monoidal closed category
(like $\boto$), see
\cite{Ben}, \cite{Bor}, or for a recent gentle introduction focusing
on applications, \cite[Section 4]{ACT}.

The opposite $P^\op$ of a poset $P$, has the same elements but order
relation reversed. The elements in the distributive lattice $\di{P}$
associated to $P$ identifies as pairs $(I,F)$, called {\em cuts},
where $I$ is a down-set in $P$
and $F$ the complement up-set. There is then a duality
between $\di{P}$ and $\di{P^\op}$ sending $(I,F)$ to $(F^\op, I^\op)$.
We call two such pairs dual or Alexander dual (as is common in
combinatorial commutative algebra).

Denote by $\Pro(P,Q)$ the profunctors $P \promap Q$. This is again a
a partially ordered set and the opposite of this poset is $\Pro(Q,P)$.
The basic notions we introduce associated to a profunctor $f : P \promap Q$
between posets are the notions of its {\it graph} $\Gamma f$ and
its {\it ascent} $\Lambda f$. These are dual notions in the sense
that if $g : Q \promap P$ is the dual profunctor, the graph of $f$ equals
the ascent of $g$. Let $UP$ and $UQ$ denote the underlying sets of $P$
and $Q$.
In $\Pro(P,Q)$ let $\cI$ be a down-set and $\cF$ its complement
up-set, so $(\cI,\cF)$ is a cut for $\Pro(P,Q)$. Let
$\cF_\Lambda$ be the {\em up-set} in the Boolean lattice of all subsets of
$UQ \times UP^\op$ generated by the ascents $\Lambda f$ for $f \in \cF$.
Let $\cI_\Gamma$ be the {\em down-set} in this Boolean lattice generated by
the complements of the graphs $\Gamma f$ for $f \in \cI$.

\begin{figure}
\begin{tikzpicture}
  {    
\draw[color=blue] (1,1) circle (3pt);
\filldraw[color=black] (2,1) circle (1pt);
\filldraw[color=black] (3,1) circle (1pt);
\filldraw[color=black] (4,1) circle (1pt);
\filldraw[color=black] (5,1) circle (1pt);

\filldraw[color=red] (1,2) circle (3pt);
\filldraw[color=red] (2,2) circle (3pt);
\draw[color=blue] (3,2) circle (3pt);
\filldraw[color=black] (4,2) circle (1pt);
\filldraw[color=black] (5,2) circle (1pt);

\filldraw[color=black] (1,3) circle (1pt);
\filldraw[color=black] (2,3) circle (1pt);
\draw[color=blue] (3,3) circle (3pt);
\filldraw[color=black] (4,3) circle (1pt);
\filldraw[color=black] (5,3) circle (1pt);

\filldraw[color=black] (1,4) circle (1pt);
\filldraw[color=black] (2,4) circle (1pt);
\filldraw[color=red] (3,4) circle (3pt);
\filldraw[color=red] (4,4) circle (3pt);
\draw[color=blue] (5,4) circle (3pt);

}
\end{tikzpicture}
\caption{Profunctor from $[5]$ to $[4]$}
\label{fig:in-pro}
\end{figure}
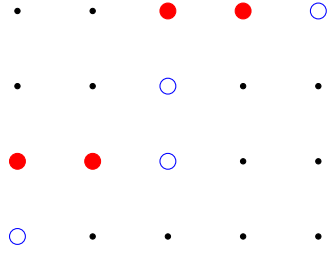

\begin{example}
  In Figure \ref{fig:in-pro} the red 
  discs give the graph of a profunctor
  $[5] \promap [4]$. The blue circles give the ascent of this profunctor.
  The graphs and the ascents are then subsets of $U[4] \times U[5]^{\op}$.
  The up-set $\cF_\Lambda$ is the up-set of the Boolean lattice
  of subsets of $U[4] \times U[5]^\op$ generated by the ascents of
  $f \in \cF$. The down-set $\cI_\Gamma$ is the down-set of
  this Boolean lattice generated by the complements of graphs of $f \in \cI$.
\end{example}
  
Our main theorem states the following.

\medskip
\noindent {\bf Theorem \ref{thm:AD-duality}.}\!\! (Preserving the cut)
  {\em Let $P$ and $Q$ be well-founded posets,
  and $(\cI,\cF)$ a cut for
  $\Pro(P,Q)$. Then $(\cI_\Gamma,\cF_\Lambda)$ is a cut
  for the Boolean lattice of subsets of $UQ \times UP^\op$.}

\medskip
\noindent{\em Example} 0.1 {\em continued.} This says
  that given any subset $S$ of $U[4] \times U[5]^\op$, then either
  $S$ contains an ascent $\Lambda f$ for $f \in \cF$, or the complement
  $S^c$ contains a graph $\Gamma f$ for $f \in \cI$. These two cases
  are also mutually exclusive.

\medskip
  \medskip
  This theorem has alternative formulations in Theorem \ref{thm:AD-ad}
  asserting that two up-sets are Alexander dual, with applications
  to Stanley-Reisner theory, and in Theorem
  \ref{thm:AD-lala} asserting that the map sending the ideal $\cI$ to
  the ideal $\cI_\Gamma$ respects
  the duality on profunctors. In Theorem \ref{thm:monwell-duality}
  we give a version
  with conditions on $P$ and $Q$ ensuring that $\Gamma f$ and $\Lambda f$
  are always finite sets, suitable for applications to monomial
  ideals, Section \ref{sec:mon}. An open problem is if a more functorial
  formulation is possible, Problem \ref{pro:AD:functorial}.

\medskip
Although we develop a general theory here, our original motivation came
from applications related to commutative algebra.

\medskip
\noindent{\em Applications to Stanley-Reisner theory.} When $P$ and $Q$
are finite posets we get general constructions, Section
\ref{sec:AD-LP}, of Alexander dual
squarefree monomial ideals, generalizing isotonian ideals and letterplace
and co-letterplace ideals, \cite{EHM},\cite{FGH}, \cite{F-LP},
\cite{Ju-Survey}, and \cite{JuKaMa}. In particular, when $Q$ is a
chain these constructions have given very large classes of simplicial
balls and spheres, \cite{DFN}. 

\medskip
\noindent{\em Applications to order preserving maps
  $f : \NN \pil \NN \cup \{\infty \}$.}
The profunctors from $\NN$ to $\NN$ identify as order preserving maps from
$\NN$ to the distributive lattice $\di{\NN}$, and the latter 
identifies as $\NN \cup \{\infty \}$. 
  Profunctors $f : \NN \promap {\NN}$ are the topic of
  many our examples. These benefit the reader with a
  quick access to many of our notions and results:
  \ref{eks:poset-NN}, \ref{eks:poset-NNY},
  \ref{eks:PQ-snake}, \ref{eks:PQ-NN}, \ref{rem:PQ-tame}, \ref{eks:AD-NN},
  \ref{eks:AD:53},
\ref{rem:top-real}, \ref{eks:monwell-fr}, \ref{eks:monwell-gr},
\ref{eks:monwell-f}.
See also the end of Section \ref{sec:mon}.

Injective order preserving maps $f : \NN \pil \NN$ form the so called
increasing monoid, which has gained recent interest.
In \cite{NaRo} Nagel and R\"omer show that
ideals in the infinite polynomial ring invariant for the
increasing monoid, have an essentially finite Gr\"obner
basis, thereby generalizing previous results for the symmetric group.
In \cite{GuSn} G\"unt\"urk\"un and Snowden studies in depth the
representation theory of the increasing monoid.
Note that the injective order preserving maps
$g : \NN \pil \NN$ are in bijection with order preserving maps
$f : \NN \pil \NN$ by $g = f + \id - 1$. Order preserving maps
$f : \NN \pil \NN$ also occur in the definition of the bicylic
semi-group \cite{EbSe}, a basic notion in inverse semi-group theory.
Our main application is in \cite{Fl-St}, relating profunctors
$f : \NN \promap \NN$ (i.e. order preserving maps
$f: \NN \pil \NN \cup \{ \infty \}$) to
the duality theory of strongly stable ideals in the the infinite
polynomial ring $k[x_i]_{i \in \NN}$.

\medskip
In order for the substantial parts of our theory, related to graphs and ascents,
to work well we must have certain conditions
on the posets $P$ and $Q$. Our weakest condition is that they
are {\em well-founded}. For our applications to
polynomial rings,
we work in the class of {\em natural} posets, Section \ref{sec:natural}.
These are posets where all
anti-chains are finite and for every $x$ in the poset the principal down-set
$\doa x$  is finite.
This is a subclass (closer to natural numbers) of
{\em well partially ordered sets}.

Another feature we introduce is a topology on $\Pro(P,Q)$, Section
\ref{sec:top}, in
particular on $\Pro(\NN, \NN)$. This
is needed for our applications to commutative algebra. For $\Pro(\NN,\NN)$
a basis for the topology consists of intervals $[f,g]$ where
i) the image $f(\NN)$ is contained in a finite interval,
and ii) $g(p) = \infty$ for all but a finite set of $p$'s.

\medskip
Organization of article:

\begin{itemize}
  \item[\ref{sec:basic}.] {\sc Preliminaries on posets.}
    Notions for posets are recalled, most significantly
    cuts for posets and the associated
distributive lattice. It relates
to simplicial complexes, Alexander duality, and Stanley-Reisner rings.
\item[\ref{sec:pro}.] {\sc Profunctors between posets.}
  We introduce these and develop basic theory.
  \item[\ref{sec:AD}.] {\sc The graph, the ascent and preserving
      the cut.} We give
    our main theorem on Alexander duality, together with variants.  
  \item[\ref{sec:AD-LP}.] {\sc Applications to finite posets and Alexander
      duality.} We connect to commutative algebra
and get Alexander dual ideals in finite dimensional polynomial rings.
\item[\ref{sec:top}.] {\sc Topology on $\Pro(P,Q)$.}
We define the topology and in particular look at {\it interior}
open down-sets.
\item[\ref{sec:natural}.] {\sc Profunctors between natural posets.}
 We consider natural posets $P$ and $Q$ and
 investigate the topology in this setting. We show an 
 open down-set is also closed (clopen down-set), if and only if it is finitely
 generated.
\item[\ref{sec:monwell}.] {\sc Natural posets and finite type cuts.}
  We give the version of the main theorem for natural posets, suitable
  to get Alexander dual ideals in (infinite dimensional) polynomial rings.
\item[\ref{sec:mon}.] {\sc Monomial ideals.}
When $Q = \NN$ we get monomial ideals
in the polynomial ring generated by $x_p, p \in P$. We briefly
indicate the applications to strongly stable ideals in \cite{Fl-St}
when $P = Q = \NN$.
\end{itemize}

\medskip
\noindent{\em Note.} We let $\NN = \{1,2,3, \ldots\}$. We only use the ordered
structure on this so we could equally well have used
$\NN_0 = \{0,1,2,3, \cdots \}$.
Only in the last Section \ref{sec:mon} do we,
in a somewhat different setting, use
the commutative monoid structure and then we explicitly write $\NN_0$. 

\medskip
\noindent {\em Acknowledgement.} We are grateful to an anonymous referee for
suggestions concerning notation and pointers to the literature.

\medskip
\noindent {\em Data availability statement.} There are now associated data
to this manuscript.

\section{Preliminaries on posets} \label{sec:basic}
We give basic notions and constructions concerning posets: down-sets, up-sets,
dualities, distributive lattices, simplicial complexes. We also
recount the algebraic notions of Stanley-Reisner ideal and ring.

\subsection{Down- and up-sets in $P$}

Let $P$ be a partially ordered set. The {\it opposite poset} $P^{\op}$
has the same underlying set as $P$ but with order relation
$\leq_\op$ where $p \leq_\op q$ if $p \geq q$ in $P$.

A {\it down-set $I$} of $P$ is a subset of $P$ closed under taking
smaller elements. An {\it up-set $F$} in $P$ is a subset of $P$ closed
under taking larger elements. If $I$ and $F$ are complements of each other,
we call $(I,F)$ a {\em cut} 
for $P$. Since each of $I$ and $F$ determine each other, we sometimes
denote this as $(-,F)$ if we focus on $F$, and similarly with $(I,-)$.
Down-sets are sometimes called {\it order ideals} and up-sets {\it order
  filters}, whence the suggestive notation $I$ and $F$. (The single term
{\it ideal} is usually reserved for order ideals in lattices closed
under joins.)

An element $p \in P$ induces the principal up-set $\upa p$ consisting
of all $p^\prime$ with $p^\prime \geq p$, and a principal down-set $\doa p$
consisting of all $p^\prime \leq p$. 

\begin{definition}
  The {\em Alexander dual} (or just {\em dual}) of the cut $(I,F)$ for $P$
  is the cut
$(F^\op, I^\op)$ for $P^\op$. The Alexander dual of the down-set $I$
is the down-set $J = F^\op$ of $P^\op$, and the Alexander dual of the
up-set $F$ is the up-set $G = I^\op$.
\end{definition}

\subsection{The distributive lattice}
If $P$ and $Q$ are two partially ordered sets, a map
$f : P \pil Q$ is {\it isotone} if it is order-preserving, i.e.
$p_1 \leq p_2$ implies $f(p_1) \leq f(p_2)$. 
We denote by $\Hom(P,Q)$ the
set of all isotone (order-preserving) maps $f: P \pil Q$. It is itself
a partially ordered set (an internal Hom)
by $f \leq g$ if $f(p) \leq g(p)$ for every $p \in P$.
The opposite poset $\Hom(P,Q)^{\op}$ naturally identifies as
$\Hom(P^{\op}, Q^{\op})$.
The category of posets forms a closed symmetric monoidal category and
so for any three posets $P,Q,R$ we have:

\begin{equation} \label{eq:PQR}
  \Hom(P \times Q, R) \iso \Hom(P,\Hom(Q,R)).
\end{equation}

\medskip

Denote by  $\boto $ the ordered set $\{ 0 < 1 \}$.
The (complete) distributive lattice associated to $P$ is
$\di{P} = \Hom(P^{\op},\boto)$.
Given an $f \in \di{P}$, the elements $p$ in $P$ such that $p^\op$ maps to
$1 \in \boto$, constitute a down-set $I$ in $P$. The complement up-set
$F$ in $P$ consists of those $p \in P$ such that $p^\op$ maps to $0$.
An element $f$ of $\di{P}$ may thus be identified either with:
\begin{itemize}
\item A down-set $I$ of $P$,
\item An up-set $F$ of $P$,
\item A cut $(I,F)$ for $P$.
\end{itemize}
We shall usually identify elements of $\di{P}$ with the down-sets $\DP$.
In categorical terms $\DP$ is the free cocompletion of $P$.
Thus an element $i \in \di{P}$ is a down-set $I$ of $P$.
However we sometimes will consider the elements of $\di{P}$ to be
cuts $(I,F)$ for $P$.
We speak
of a cut $(I,F)$ {\it in} $\di{P}$,
or equivalently a cut $(I,F)$ {\it for} $P$.
The cuts for $P$ are then  ordered by
\[ (I,F) \leq (J,G) \text{ if and only if } I \sus J \text{ or equivalently }
  F \supseteq G. \]

The distributive lattice
$\di{P}$ has a unique maximal element, denoted $\ifst$. It sends every
$p^{\op}$ to $1$, and corresponds to the cut $(I,F) = (P,\emptyset)$.
\begin{example} \hskip 2mm \hbox{}  \label{eks:poset-NN}

  \begin{itemize}
  \item $P = [n] = \{ 1 < 2 <  \cdots < n\}$ has cuts
    $(I,F)$ where $I = \{1,2,\ldots, i-1\}$ and $F = \{i, \ldots, n\}$
    for $i = 1, \ldots, n+1$. Thus $\di{P} \iso
    [n+1] = [n] \cup \{ \infty \}$ where $\infty = n+1$.
  \item Any set $A$ may be considered a discrete poset (only relations
    are $a \leq a$). Then $\di{A}$ is the Boolean lattice on $A$. It
    consists of subsets $S \sus A$. We identify such a subset with
    the cut $(S,S^c)$. (So for instance the cut $(-,T)$ identifies
    as the subset $T^c$ in $A$.)
  \item If $P = \NN$ the natural numbers, then $\di{\NN} = \NN \cup
    \{ \infty \}$.
  \end{itemize}
\end{example}

Given an element $p$ in $P$ we get a map
\[ \hat{p} : P^\op \pil \boto, \]
where if $p^\prime \geq p$ we send $p^{\prime \op} \mapsto 0$ and all other
elements of $P^\op$ to $1$.
It corresponds to the cut $((\upa p)^c,\upa p)$ for $P$ (where $()^c$ denotes
the complement set).
\begin{figure}
    \begin{tikzpicture}
 \fill[blue!5] (1,1) ellipse (1 and 2.5);     
\draw (1,1)--(1.8,2.5);
\draw (1,1)--(0.2,2.5);
\draw node[anchor=west] at (1,1) {$p^{}$};
\draw node[anchor=west] at (0.7,2.3) {$p^{}$};
\draw node at (1.1,1.7) {$0$};
\draw node at (0.8,0.5) {$1$};
\filldraw (1,1) circle (2pt);
\filldraw (0.7,2.3) circle (2pt);
\end{tikzpicture}
\caption{The poset $P$ and the cut $\di{P}$}
\label{fig:settPoset}
\end{figure}
This gives a distinguished injective poset map
\begin{equation}  \label{eq:pos-Yvar}
  P  \lpil  \di{P} = \Hom(P^{\op}, \boto), \quad
  p  \mapsto \hat{p},  \text{ corresponding to the cut } ((\upa p)^c,\upa p)
\end{equation}
see Figure \ref{fig:settPoset}.
Note that the image of the map \eqref{eq:pos-Yvar} is in
$\di{P} \backslash \{ \ifst \}$.
  
\begin{remark} A poset may be considered a $\boto$-category for
  the symmetric monoidal closed category $\boto$.
  Where
  \[ \Hom_\boto(p^\prime,p) = \begin{cases} 1 &  p^\prime \leq p \\
      0 & p^\prime \not \leq p
      \end{cases}. \]
  The map \eqref{eq:pos-Yvar} is not the Yoneda embedding
  \begin{equation} \label{eq:pos-Ycov}
    P \pil \Hom(P^\op , \boto), \quad p \mapsto  \Hom_\boto(-,p).
    \end{equation}
  Rather the map \eqref{eq:pos-Yvar} is derived as follows.
  One has the Yoneda embedding:
   \[ P^\op \pil \Hom(P , \boto), \quad p \mapsto  \Hom_\boto(p,-). \]
   Taking the opposite of this we get:
   \begin{equation} \label{eq:pos-Ycon}  P \pil \Hom(P,\boto)^\op.
   \end{equation}
   This is the {\it co-Yoneda} embedding, \cite{Ba}.
   Note that $\boto^\op \iso \boto$ by sending $0^\op \mapsto 1$ and
   $1^\op \mapsto 0$. So
   \begin{equation} \label{eq:pos-boto}
     \Hom(P,\boto)^\op =   \Hom(P^\op, \boto^\op)  \iso \Hom(P^\op, \boto)
     = \di{P}. \end{equation}
   In general however,
   for instance ordinary categories with ${\bf set}$
   instead of $\boto$, the Yoneda and co-Yoneda embedding map to different
   categories. 
   Composing \eqref{eq:pos-Ycon} above with \eqref{eq:pos-boto} 
   we get the embedding \eqref{eq:pos-Yvar}.
 \end{remark}

 \begin{example} \label{eks:poset-NNY}
   Let $P = [n]$. The map \eqref{eq:pos-Yvar}, the co-Yoneda embedding, is 
   \[ [n] \lpil \diw{[n]} = [n+1], \quad i \mapsto i. \]
   On the other hand the Yoneda embedding \eqref{eq:pos-Ycov}
   \[ [n] \lpil \diw{[n]} = [n+1], \quad i \pil i+1.  \]
   In particular the top element $n \mapsto n+1 = \infty$.

   In more abstract terms, the Yoneda embedding embeds a chain of length
   $n$ into the upper part of a chain of length $n+1$, while the co-Yoneda
   embedding embeds it into the lower part. 
 \end{example}


\subsection{Down- and up-sets in $\di{P}$}
Let $\cI $ be a down-set for $\di{P}$, and $\cF $ the complement up-set
of $\di{P}$. So $\cI $ consists of cuts $(I,F)$ closed under forming cuts with
smaller $I$'s, and $\cF $ consists of cuts $(I,F)$ closed under forming
cuts with larger $I$'s (or equivalently smaller $F$'s).
Then $(\cI,\cF)$ is a cut {\em for} $\di{P}$ (note again terminology:
$(I,F)$ is a cut {\em in} $\di{P}$). Also $(\cI,\cF)$ is a cut {\em in}
$\ddhat{P}$.

\begin{lemma} \label{lem:posP-IFd}
  A cut $(I^\prime, F^\prime)$ in $\di{P}$ is in $\cI$ iff
  $F^\prime \cap I \neq \emptyset$ for every $(I,F)$ in $\cF$.
\end{lemma}

\begin{proof}
That $(I^\prime,F^\prime)$ is in $\cI$ means that we cannot find any
  $(I,F)$ in $\cF$ such that $I \sus I^\prime$. Alternatively
  $F^\prime \cap I \neq \emptyset $ for each $(I,F)$ in $\cF$.
\end{proof}

For later use we look closer at Alexander duality for $\di{P}$
and $\ddhat{P}$. A cut $(I,F)$ {\it for } ${P}$, i.e. a cut in $\di{P}$,
gives a dual cut $(F^{\op},I^\op)$ for ${P}^\op $.
A cut $(\cI,\cF)$ {\it for}
$\di{P}$, i.e. a cut in $\ddhat{P}$ gives a dual cut $(\cF^\op, \cI^\op)$
for $(\di{P})^\op = P^{\, \widehat{\op}}$, so this is a cut in $P^{\, \ddbhat{\op}}$.
For $u = (I,F) \in \cI \sus \di{P}$, then
$u^\op \in \cI^\op$ is
\[ u^\op = (I,F)^\op = (F^\op, I^\op) \in P^{\, \widehat{\op}}. \]
So when we take the dual of the cut $(\cI,\cF)$ for $\di{P}$,
we not only get a switch $(\cF^\op, \cI^\op)$ but also the elements
$u = (I,F)$ of $\cI$ or $\cF$ are switched, to $u^\op = (F^\op, I^\op)$.

 \subsection{Finite type cuts}

 The elements of the distributive lattice $\di{P}$ identify as cuts $(I,F)$.

 \begin{definition} {\hskip 2mm}
   
   \begin{itemize}
\item $\difi{P}$ is the sublattice of $\di{P}$
consisting finite cuts: cuts $(I,F)$ where $I$ is finite.
\item $\diff{P}$ is the sublattice of $\di{P}$ consisting of
cofinite cuts: cuts $(I,F)$ where $F$ is finite.
\end{itemize}
\end{definition}

When $P$ is an infinite discrete poset, an infinite set, such cuts $(I,F)$
of $P$ are a standard example of infinite Boolean algebras,
called {\em finite-cofinite} algebras.

\begin{definition} \label{def:sett-AD}
  A {\it finite type} $\di{P}$-cut
  is a pair $(\cI,\cF)$ where
  $\cF$ is an up-set for $\difi{P}$ and $\cI$ an ideal
  for $\diff{P}$ such that the following holds.
  \begin{itemize}
    \item[1.] $ (I,F) \in \cI \text{ if and only if } F \cap J \neq \emptyset
        \text{ for every } (J,G) \in \cF. $
      \item[2.] $ (J,G) \in \cF \text{ if and only if } F \cap J \neq \emptyset
        \text{ for every } (I,F) \in \cI. $
    \end{itemize}
  \end{definition}

\ignore{\begin{definition} \label{def:sett-AD}
  Let $\cF$ be a up-set for $\dif{P}$ and $\cG$ a up-set
  for $\dif{P}^\op$. These up-sets are {\it Alexander dual} if:
  \begin{itemize}
    \item[1.] $(J,G) \in \cG \text{ if and only if } J^\op \cap I \neq \emptyset
        \text{ for every } (I,F) \in \cF .$
      \item[2.] $ (I,F) \in \cF \text{ if and only if } I^\op \cap J \neq \emptyset
        \text{ for every } (J,G) \in \cG .$
    \end{itemize}
  \end{definition}
}
  Note that for $P$ finite
  then 1 and 2 are equivalent. If $P$ is infinite, one of the above will
  in general not imply the other. The point of having both fulfilled
  is that $\cF$ determines $\cI$ by 1, and vice versa by 2. If only
  1 holds then $\cF$ determines $\cI$, but one may not be able to reconstruct
  $\cF$ from $\cI$. 

  In Section \ref{sec:monwell} we construct such finite type cuts for infinite
  posets.



 \subsection{Simplicial complexes and Stanley-Reisner rings}

 Let $A$ be a set. A simplicial complex $X$ on $A$ is a set of subsets of
 $A$ closed under taking smaller subsets, i.e. if $I \in X$ and $J \sus I$,
 then $J \in X$.
 
 The set $A$ may be considered as a poset with the discrete poset structure,
 i.e. the only comparable elements are $a \leq a$ for $a \in A$. 
 Then $\di{A}$ identifies as the Boolean lattice on $A$ (see Example
 \ref{eks:poset-NN}),
 consisting of all subsets of $A$.
  A cut $(\cI,\cF)$ for
  $\di{A}$ corresponds precisely to a simplicial complex $X$.
  The elements $I$ in $X$ give the cuts $(I,I^c)$
  in $\cI$.

  \medskip
  The Alexander dual simplicial complex $Y$ of $X$
  consists of all the complements $I^c$ of subsets $I \sus A$ such that
  $I$ is not in $X$. 
The Alexander dual cut $(\cF^\op, \cI^\op)$ for $\di{A}^{\op} \iso \di{A}$ then
  corresponds to $Y$: The cuts $(F^\op, I^\op)$ in $\cF^\op$ give precisely
  the elements $F^\op = (I^\op)^c$ in $Y$.

  \medskip
  Denote by $k[x_A]$ the polynomial ring in the variables $x_a$ for $a \in A$.
  When $A$ is finite, to the simplicial complex $X$ corresponding to the cut
  $(\cI,\cF)$, we associate a monomial ideal $I_X$ in $k[x_A]$,
  the Stanley-Reisner ideal of $X$. It is generated by monomials
  $x_I = \prod_{i \in I}x_i$ for $(I,F) \in \cF$. These are the subset $I$ of
  $A$ such that $I$ is {\em not} in the simplicial complex $X$.
The monomials in the Alexander dual Stanley-Reisner ideal $I_Y$
  are then precisely those monomials which have non-trivial common divisor
  with every monomial in $I_X$, by the characterization of Lemma
  \ref{lem:posP-IFd}.

\medskip
  When $A$ is infinite, we still have a polynomial ring $k[x_A]$. But
  the construction above does not give meaning if $\cF$ is non-empty, since
  there would be $(I,F)$ in $\cF$ with $I$ infinite.
  However if $(\cI,\cF)$ is a finite type $\di{A}-$cut,
  it gives rise to a Stanley-Reisner ideal $I_X$ of $k[x_A]$, as above.

\section{Profunctors between posets} \label{sec:pro}
We introduce profunctors $P \promap Q$ between posets. Such a profunctor has
a dual $Q \promap P$ and we investigate this correspondence.
 For an introduction to profunctors, see \cite[Chap.4]{ACT}.
 See also \cite[Section 7]{Bor} and \cite{Ben} where they are called
 distributors.

\subsection{Duality on profunctors}
A {\it profunctor} $P \promap Q$
is simply a poset homomorphism
$P \pil \di{Q}$.
By the adjunction
\begin{eqnarray*}  \Hom(P,\di{Q}) & = &\Hom(P, \Hom(Q^\op, \boto))\\
                                 & = & \Hom(P \times Q^\op, \boto) \\
                                 & = & \Hom((Q \times P^\op)^\op, \boto)
                                       = (Q \times P^\op)^\dual.
\end{eqnarray*}
Thus a profunctor is equivalently an isotone map $P \times Q^\op \pil \boto$
and this is often taken as the definition. It is also equivalently an
element of the distributive lattice $(Q \times P^\op)^\dual$, and so
corresponds to a cut $(I,F)$ for $Q \times P^\op$.
(Our convention differs somewhat from \cite{ACT}, since there a profunctor
$P \promap Q$ 
corresponds to an isotone map $P^\op \times Q \pil \boto$.)

In particular if $Q = B$ and $P = A$ are sets, this is simply a subset
of $B \times A^\op$ or a relation between the sets $A$ and $B$.
Profunctors are also called distributors or bimodules in the literature.
We denote the set of profunctors $P \promap Q$ as
$\Pro(P,Q) (= \Hom(P,\di{Q}))$. It is itself a partially ordered set,
in fact a distributive lattice. 

\begin{remark} Profunctors can be taken to be the morphisms between
  two ordered sets. They may be composed and form a category.
  R.Rosebrugh and R.J.Wood show in \cite{CCD4} that
  this category is equivalent to
  the category of {\it totally algebraic lattices}, those of the form
   $\di{P}$ for some poset $P$,  with supremum-preserving isotone
   maps. (There the term {\it ideal} is used for profunctor.)

   Order theory may also be done in a more categorical setting, for objects
   in a topos as in \cite{CCD4}, or even more general categories \cite{Ca}.
   The result of Rosebrugh and Wood mentioned above also has a more
   general categorical formulation \cite{RW04} in terms of a monad on a
   category where idempotents split.
\end{remark}

When $(I,F)$ is the cut of $Q \times P^\op$ corresponding to $f$,
the down-set $I$ may be considered as a relation defining the profunctor.
We then write $qfp$ when $(q,p^\op) \in I$.

\begin{lemma} \label{lem:posPQ-qfp}
  Given a profunctor $f : P \promap Q$.
  \begin{itemize}
  \item[a.] The following are equivalent:
    i) $q \in f(p)$ and ii) $qfp$.
    \item[b.] The following are equivalent:
      i) $q  \in f(p)^c$, ii) $f(p) \leq \hat{q}$, and iii) $\neg \, qfp$.
      \end{itemize}
\end{lemma}

\begin{proof} Let $f$ correspond to $P \pil \di{Q} = \Hom(Q^\op, \boto)$. 
  That $q \in f(p)$, the latter a down-set,  means that $q^\op  \mapsto 1$. This
  says $(p,q^\op) \mapsto 1$, and so $(q,p^\op)$ is in the
  down-set $I$ corresponding to $f$.

  The element $\hat{q}$ corresponds to the cut $((\upa q)^c, \upa q)$.
  That $f(p)\leq  \di{q}$ then means that $f(p) \sus (\upa q)^c$, or
  equivalently $q \not \in f(p)$.
\end{proof}

\medskip 
Since $\Pro(P,Q)$ identifies as $(Q \times P^\op)^\dual$,
the following is seen to be natural, by taking opposites.

\begin{lemma} \label{lem:settD} Let $P,Q$ be posets.
There is a a natural isomorphism of
  posets
  \[ \Pro(P,Q)^{\op} \overset{D}\iso \Pro(Q,P). \]

\end{lemma}

\begin{proof}
 First
  \[ \Hom(P,\di{Q})^{\op} \iso \Hom(P^\op, (\di{Q})^{\op}) \iso
    \Hom(P^{\op}, \Hom(Q,\boto^\op)). \]
Using that $\boto^{\op}$ naturally is isomorphic
  to $\boto$, and the adjunction \eqref{eq:PQR} this further becomes
  \begin{align*}
    \Hom(P^{\op}, \Hom(Q,\boto)) \iso \Hom(P^{\op} \times Q, \boto)
   & \iso \Hom(Q,\Hom(P^{\op}, \boto)) \\
   & \iso \Hom(Q,\di{P}).
    \end{align*}
\end{proof}

Here is more detail on the duality $D$.

\begin{lemma} \label{lem:posPQ-Dfq}
  Given a profunctor  $f : P \promap Q$ and its dual $g = Df : Q \promap P$.
  \begin{itemize}
  \item[a.] $q f p$ if and only if $\neg \,  pgq $,
  \item[b.] $f(p) = \{ q \in Q \, | \, g(q) \leq \di{p}\}$.
  \end{itemize}
    In particular
      $f(p) = \ifst$ if and only if $g(q) \leq \di{p}$ for every $q \in Q$.
 \end{lemma}

 \begin{proof}
   If $f$ corresponds to the cut $(I,F)$ of $Q \times P^\op$, the
   dual profunctor $g$ corresponds to the cut $(F^\op, I^\op)$ of
   $P \times Q^\op$. The statements in a are equivalent to $(q,p^\op) \in I$.

   For part b, the condition $q \in f(p)$ is equivalent to $g(q) \leq \di{p}$
   by part a and Lemma \ref{lem:posPQ-qfp}.
   \end{proof}

\begin{example} \label{eks:PQ-snake}
  Let $P = Q = \NN = \{1,2,3,\cdots \}$ and consider a profunctor
  $f : \NN \promap {\NN}$ (see Example \ref{eks:poset-NN}) which
  is an isotone map $f : \NN \pil \NN \cup \{\infty \}$. Let its values be
  \[ 2,2,4,5,5,7, \cdots . \]
In Figure \ref{fig:settfDf} the graph of $f$ are marked with red discs
  \begin{tikzpicture}
    \filldraw[color=red] (1,2) circle (2pt);
    \end{tikzpicture}.
    We fill in with blue circles
   \begin{tikzpicture}
    \draw[color=blue] (1,2) circle (2pt);
  \end{tikzpicture}
  to make a connected ``snake'', starting at position $(1,1)$. 
  The graph of the dual map $g = Df$ is given by the blue circles by considering
  the vertical axis as the
argument for $g$. The values of $g$ are
\[ 1,3,3,4,6,6, \cdots .\]
Observe that for a profunctor $f : \NN \promap \NN$ (which identifies
as an isotone $f : \NN \pil \di{\NN}$), then $f(1) \geq 2$ iff
the dual map $g = Df$ has $g(1) = 1$. Hence there are no self-dual maps
$f$. 

\begin{figure}
\begin{tikzpicture}
\draw (1,1)--(2,1)--(3,1)--(4,1) --(5,1)--(6,1)--(7,1);
\draw (1,1)--(1,2)--(1,3)--(1,4) --(1,5)--(1,6)--(1,7);
\draw node[anchor=north] at (1,1) {1}
node[anchor=north] at (2,1) {2}
node[anchor=north] at (3,1) {3}
node[anchor=north] at (4,1) {4}
node[anchor=north] at (5,1) {5}
node[anchor=north] at (6,1) {6};

\foreach \x  in {1,...,6}{
        \draw (\x cm, 1cm + 1pt) -- (\x cm,1cm-1pt);
};
        
\draw node[anchor=east] at (1,1) {1}
node[anchor=east] at (1,2) {2}
node[anchor=east] at (1,3) {3}
node[anchor=east] at (1,4) {4}
node[anchor=east] at (1,5) {5}
node[anchor=east] at (1,6) {6};

\foreach \x  in {1,...,6}{
        \draw (1cm + 1pt,\x cm) -- (1cm-1pt,\x cm);
      };
      
\filldraw[color=red] (1,2) circle (2pt);
\filldraw[color=red] (2,2) circle (2pt);
\filldraw[color=red] (3,4) circle (2pt);
\filldraw[color=red] (4,5) circle (2pt);
\filldraw[color=red] (5,5) circle (2pt);
\filldraw[color=red] (6,7) circle (2pt);

\draw[color=blue] (1,1) circle (2pt);
\draw[color=blue] (3,2) circle (2pt);
\draw[color=blue] (3,3) circle (2pt);
\draw[color=blue] (4,4) circle (2pt);
\draw[color=blue] (6,5) circle (2pt);
\draw[color=blue] (6,6) circle (2pt);
\end{tikzpicture}
\caption{}
\label{fig:settfDf}
\end{figure}

The profunctor $f$ corresponds to the cut $(I,F)$ for
$\NN \times \NN^\op$ (where $\NN$ corresponds to the $y$-axis and $\NN^\op$
to the (reversed) $x$-axis)
where the up-set $F$ is given by filling in red discs vertically above those
in the graph, see Figure \ref{fig:settIF}, 
and the ideal $I$ is given by filling in blue circles
to the right of those which are present in Figure \ref{fig:settfDf}.

\begin{figure}
\begin{tikzpicture}
\draw (1,1)--(2,1)--(3,1)--(4,1) --(5,1)--(6,1)--(7,1);
\draw (1,1)--(1,2)--(1,3)--(1,4) --(1,5)--(1,6)--(1,7);
\draw node[anchor=north] at (1,1) {1}
node[anchor=north] at (2,1) {2}
node[anchor=north] at (3,1) {3}
node[anchor=north] at (4,1) {4}
node[anchor=north] at (5,1) {5}
node[anchor=north] at (6,1) {6};

\foreach \x  in {1,...,6}{
        \draw (\x cm, 1cm + 1pt) -- (\x cm,1cm-1pt);
};
        
\draw node[anchor=east] at (1,1) {1}
node[anchor=east] at (1,2) {2}
node[anchor=east] at (1,3) {3}
node[anchor=east] at (1,4) {4}
node[anchor=east] at (1,5) {5}
node[anchor=east] at (1,6) {6};

\foreach \x  in {1,...,6}{
        \draw (1cm + 1pt,\x cm) -- (1cm-1pt,\x cm);
      };
      
\filldraw[color=red] (1,2) circle (2pt);
\filldraw[color=red] (2,2) circle (2pt);
\filldraw[color=red] (1,3) circle (1.3pt);
\filldraw[color=red] (2,3) circle (1.3pt);
\filldraw[color=red] (2,4) circle (1.5pt);
\draw node[color=red] at (2,5) {\vdots};

\filldraw[color=red] (3,4) circle (2pt);
\filldraw[color=red] (3,5) circle (1.3pt);

\filldraw[color=red] (4,5) circle (2pt);
\filldraw[color=red] (4,6) circle (1.3pt);

\filldraw[color=red] (5,5) circle (2pt);
\filldraw[color=red] (6,7) circle (2pt);
\filldraw[color=red] (5,6) circle (1.3pt);
\filldraw[color=red] (5,7) circle (1.3pt);

\draw[color=blue] (1,1) circle (2pt);
\draw[color=blue] (2,1) circle (1.3pt);
\draw[color=blue] (3,1) circle (1.3pt);
\draw[color=blue] (4,1) circle (1.3pt);

\draw[color=blue] (3,2) circle (2pt);
\draw[color=blue] (3,3) circle (2pt);
\draw[color=blue] (4,2) circle (1.3pt);
\draw node[color=blue] at (5,2) {$\cdots$};

\draw[color=blue] (4,4) circle (2pt);
\draw[color=blue] (4,3) circle (1.3pt);
\draw[color=blue] (5,3) circle (1.3pt);
\draw[color=blue] (5,4) circle (1.3pt);
\draw[color=blue] (6,4) circle (1.3pt);

\draw[color=blue] (6,5) circle (2pt);
\draw[color=blue] (6,6) circle (2pt);
\end{tikzpicture}
\caption{}
\label{fig:settIF}
\end{figure}

\end{example}


\ignore{
\begin{lemma}
  Let $(\cI,\cF)$ be a cut for $\Hom(P,\di{Q})$. It gives a dual cut
  $(\cJ,\cG) = (D\cF, D\cI)$ for $\Hom(Q,\di{P})$. Then $g$ is in the dual poset
  ideal $\cJ = D\cF$ if and only if for every $f$ in $\cI$ there are $p \in P$ and
  $q \in Q$ such that
  \begin{equation} \label{eq:settfgpq}
    i. g(q) \leq \hat{p}, \quad ii. f(p) \leq \hat{q}.
    \end{equation}
\end{lemma}
  

\begin{proof}
  The map $g$ is in $D\cF$ if and only if $Dg$ is in $\cF$ and this holds iff
  there is not $f$ in $\cI$ with $g \leq f$. This means that 
  for every $f$ in $\cI$ there is a $p$ with $Dg(p)$ not $\leq f(p)$.
  The again  means that there is $q$ in $Dg(p) \backslash f(p)$. This means:
  \begin{itemize}
  \item[i.] $Dg(p)$ sends $q^{\op}$ to $u$,
  \item[ii.] $f(p)$ sending $q^{\op}$ do $d$.
  \end{itemize}
  By Lemma \ref{lem:posPQ-Dfq},
  i above is equivalent to i in \eqref{eq:settfgpq}.
\end{proof}

When $Q$ is totally ordered, then $\di{Q}$ may be considered as the union
$Q \cup \{ \infty \}$. We may then specialize the above to the following.

\begin{corollary}
  If $Q$ is totally ordered, then $g$ is in $\cJ = D\cF$ if and only if it
  is so ``small'' such that for every
  $f$ in $\cI$ there is a $p \in P$ with $f(p) \neq \infty$ and
  $g(f(p)) \leq \hat{p}$.
\end{corollary}
}

\subsection{Alexander duality}

\begin{definition}
  Recall that a profunctor $f : P \promap Q$ is an isotone
  $f : P \pil \di{Q}$. 
  \begin{itemize}
    \item For fixed $q \in Q$, the down-set 
      $\Pro_{\leq \di{q}}(P,Q)$ is the set of profunctors $f$
      such that $f(p) \leq   \di{q}$ for all $p$.
 \item For fixed $p \in P$, the down-set $\Pro^{\im(p) < \infty}(P,Q)$
  is the set of profunctors $f$ such that $f(p) < \infty$.
  \item The down-set $\Pro^{< \infty}(P,Q)$ is the set of profunctors
  $f$ such that $f(p) < \infty$ for every $p \in P$.
\end{itemize}
\end{definition}

\begin{lemma} \label{lem:PQ-ad}
  The down-sets $\Pro_{\leq \di{q}}(P,Q)$ and $\Pro^{\im(q) < \infty}(Q,P)$
  are Alexander duals.
\end{lemma}

\begin{proof}
  That $f(p) \leq \di{q}$ is equivalent to, letting $g = Df$, that
  $p \in g(q)$. That this holds for all $p \in P$ is equivalent
  to $g(q) = \infty$. The Alexander dual down-set of
  $\Pro_{\leq \di{q}}(P,Q)$ is then those maps $g$ {\em not}
fulfilling $g(q) = \infty$. 
\end{proof}

\begin{corollary} \label{cor:PQ-ad}
  If $Q$ has a maximal element, then
  $\Pro^{< \infty}(P,Q)$ and $\Pro^{< \infty}(Q,P)$ are Alexander
  dual down-sets.

  Note: By symmetry of the conclusion, this holds under the weaker assumption
  that $P$ or $Q$ has a maximal element.
\end{corollary}

\begin{proof} If $q$ is the maximal element of $Q$
  the down-sets of Lemma \ref{lem:PQ-ad} are respectively
  $\Pro^{< \infty}(P,Q)$ and $\Pro^{< \infty}(Q,P)$.
\end{proof}

\subsection{Profiles and co-profiles}

\begin{definition}
  The {\it profile} of a profunctor $f : P \pil Q$ is the cut
  $(I,F)$ for $P$ where the {\it profile up-set} $F$ consists of
  all $p \in P$ such that $f(p) = \infty$, and the {\it profile down-set}
  $I = F^c$ is the complement down-set.

The {\it co-profile} of $f$ is the cut $(J,G)$ for $Q$ where
$J$ is the union of all $f(p)$ (considered as a down-set of $Q$)
for $p \in P$, and
$G$ is the complement of $J$. 
\end{definition}


The following is immediate.

\begin{lemma} \label{lem:posPQ-profile}
  Let $f : P \promap Q$ be a profunctor.
  \begin{itemize}
\item[a.] Its profile filter
  $F = \{ p \in P \, | \, qfp \text{ for every } q \in Q \}$.
\item[b.] Its co-profile filter
  $G = \{ q \in Q \, | \, \neg pfq \text{ for every } p \in P \}$.
\end{itemize}
By Lemma \ref{lem:posPQ-Dfq}, if $g = Df$ is the dual, the co-profile of $f$
equals the profile of $g$.
\end{lemma}

We identify the following subsets of $\Pro(P,Q)$.

\begin{itemize}
\item $\Proi(P,Q)$ are the profunctors $f$ whose profile down-set $I$ is finite.
  These maps are called {\em large}. Then $f(p) = \infty$ for all but a finite
  number of $p$'s.
\item $\Prob(P,Q)$ are the $f$ whose coprofile down-set $J$ is finite.
  These maps are called {\em small}. Then there is a finite down-set
  bounding the $f(p)$, i.e. all $f(p) \sus J$. 
\item $\Prou(P,Q)$ are the $f$ which are in neither the above,
  so both the profile down-set $I$ and co-profile down-set $J$ are infinite.
\end{itemize}

  A consequence of Lemma \ref{lem:posPQ-profile} is the following.
  \begin{lemma} \hskip 1mm
  \begin{itemize}
\item[a.] The duality $D$ switches $\Proi(P,Q)$ and $\Prob(Q, P)$ and
maps $\Prou(P,Q)$ to $\Prou(Q,P)$. 
\item[b.] If $P$ is finite, then $\Proi(P,Q) = \Pro(P,Q)$.
\item[c.]  If $Q$ is finite, then $\Prob(P,Q) = \Pro(P,Q)$.
\end{itemize}
\end{lemma}

\begin{example} \label{eks:PQ-NN}
  Consider profunctors $f : \NN \promap \NN$. By Example
  \ref{eks:poset-NN} recall $\di{\NN} = \NN \cup \{ \infty \}$.
  Such a map is large
  if $f(n) = \infty$ for some $n$. It is small if $f$ is eventually constant,
  so $f(n) = c$ for all $n \geq n_0$. The maps in $\Pro^u(\NN, \NN)$
  are the maps $f : \NN \pil \NN$ which are not bounded,
  so $\lim_{n \pil \infty} f(n) = \infty$. We see the naturalness of considering
  $\Pro(\NN, \NN)$ instead of $\Hom(\NN, \NN)$: The latter is not
  self-dual while the former is. $\Pro(\NN, \NN)$ has two
  countable dual ``shores'' $\Proi(\NN, \NN)$ and $\Prob(\NN,\NN)$
  and between them an uncountable self-dual ``ocean'' $\Prou(\NN,\NN)$.
  \end{example}

  \begin{remark} \label{rem:PQ-tame}
    The small profunctors $\Prob(\NN,\NN)$, which are simply
    bounded maps $f : \NN \pil \NN$
    are in bijection with the
    tame (small) increasing monoid of \cite{GuSn} by sending the bounded
    map $f$ to $f + \id - 1$.
    \end{remark}


\subsection{Adjunctions} \label{subsec:PQ-adj}

Given an isotone map $f : P \pil Q$ it induces
a pull-back map
\[f^* : \di{Q} \pil \di{P}, \quad (J,G) \mapsto (f^{-1}(J), f^{-1}(G)). \]
This map has a left adjoint
\begin{equation*}
  f^! : \di{P}  \pil \di{Q},  \quad (I,F) \mapsto (f(I)^\downarrow, -)
\end{equation*}
where $f(I)^\downarrow $ is the smallest down-set in $Q$ containing $f(I)$.
There is also a right adjoint of $f^*$:
\begin{equation*}
  f^\ii : \di{P}  \pil \di{Q},  \quad (I,F) \mapsto (-,f(F)^\uparrow)
\end{equation*}
where $f(F)^\uparrow$ is the smallest up-set in $Q$ containing
$f(F)$.
All these maps are functorial in $P$ and $Q$.

\subsubsection{A variation on the setting}
There is a variation for these maps as follows.
There is a forgetful functor $U : \Poset \pil \Set$ by mapping a poset
to the underlying set. Composing with the natural inclusion $\Set \pil \Poset$
we get $U : \Poset \pil \Poset$. The inclusion
$i : UP \pil P$ induces the dual map $i^* : \di{P} \pil (UP)^\dual$.

Suppose we have an isotone map 
of posets
\[ g : UP \pil Q.\] (This
instead of an isotone map $f : P \pil Q$.
Note also that $g$ is really just a map of sets $UP \pil UQ$.)
We then get composites
\begin{equation} \label{eq:PQ-gu}
  g_U^\ii \, \, (\text{resp. \,} g_U^!) : \di{P} \mto{i^*}
  (UP)^\dual \xrightarrow{g^\ii \, (\text{resp. \,} g^!) } \di{Q}.
\end{equation}
Here $g_U^\ii$ sends a cut $(I,F)$ of $P$ to the cut $(J,G)$ of $Q$
where $G = g(F)^\uparrow$ is
the up-set of $Q$ generated by the $g(p)$ for $p \in F$.
Note that the above composite
$g^\ii \circ i^*$ has a left adjoint $i^! \circ g^*$ (although it will not
play a role for us).

Also note the following
\begin{equation} \label{eq:gud}
  (g_U^{!})^\op = (g^{\op}_U)^\ii : \di{P^\op} \pil \di{Q^\op}.
\end{equation}

\section{The graph, the ascent,  and preserving the cut}
\label{sec:AD}
We define the two significant notions of this article, the graph and
ascent of a profunctor $P \promap Q$, or equivalently the right
and left boundaries of the cut $(I,F)$ for $Q \times P^\op$
corresponding to this profunctor.
Then we state several versions of the main theorem of this article,
Theorem \ref{thm:AD-duality}, on preserving cuts.


\subsection{The graph and ascent}
\begin{definition}
  Denote by $\bfnu$ the minimal element in $\Pro(P,Q)$. It
  sends every $p \mapsto \emptyset$, corresponding to the cut $(\emptyset, Q)$.
Denote by $\bfen$ the maximal element in $\Pro(P,Q)$. It
  sends every $p \mapsto Q$, corresponding to the cut $(Q,\emptyset)$.
\end{definition}

Recall that for a poset $P$  then $UP$ denotes the underlying set,
considered as a discrete poset.

\begin{definition}
Let $f : P \promap Q$ be a profunctor. Its {\em ascent} is
\begin{equation} \label{eq:AD-LaDef}
  \Lambda f = \{ (q,p) \, | \, q \in f(p) \text{ but }
  q \not \in f(p^\prime) \text{ for } p^\prime < p\} \sus UQ \times UP^\op.
\end{equation}
Its {\it graph} is
\begin{equation} \label{eq:AD-GaDef} \Gamma f = \{ (q,p)\, | \, q
  \text{ minimal in the complement } 
  f(p)^c\} \sus UQ \times UP^\op.
  \end{equation}
\end{definition}

\begin{example} \label{eks:AD-NN}
  For a profunctor $f : \NN \promap \NN$, 
  see Figure \ref{fig:settfDf}, then $\Gamma f$ is the red
  discs, and $\Lambda f$ is the blue circles.
 \end{example}

Note that $\Lambda f$ and $\Gamma f$ are disjoint.
Also note i) $\Lambda \bfnu = \emptyset$ and $\Gamma \bfnu = \min Q \times UP$
where $\min$ denotes the minimal elements, ii) $\Lambda \bfen = UQ \times
\min P$ and $\Gamma \bfen = \emptyset$. 



\begin{remark} Their union $Bf = \Lambda f \cup \Gamma f$, might be called the
{\em boundary} of the cut $(I,F)$. This union will not play a role
here, but it does in \cite[Section 2]{DFN}.
\end{remark}

\subsection{Left and right boundaries}
We have seen that $\Pro(P,Q)$, the profunctors from $P$ to $Q$, identify
as $(Q \times P^\op)^\dual$. So a profunctor
$f : P \promap Q$ corresponds to a cut $(I,F)$ for $Q \times P^\op$ where 
\begin{equation} \label{eq:alexd-IF}
I = \{ (q,p) \, | \, q \in f(p) \}, \quad
F = \{ (q,p) \, | \, q \not \in f(p) \}.
\end{equation}

\begin{definition} If $(I,F)$ is a cut for $Q \times P^\op$ its
{\it left and right boundaries} are respectively 
\begin{eqnarray*}
\Lambda I & = \{ (q,p) \in I \, | \, (q, p^\prime) \not \in I \text{ for }
p^\prime < p \} \sus UQ \times UP^\op \\
\Gamma F & = \{ (q,p) \in F \, | \, (q^\prime,p) \not \in F 
           \text{ for } q^\prime < q \} \sus UQ \times UP^\op
           \end{eqnarray*}
\end{definition}

We see immediately by \eqref{eq:alexd-IF}
that if $f$ corresponds to $(I,F)$ then
$\Lambda f = \Lambda I$ and $\Gamma f = \Gamma F$.

\begin{corollary} \label{cor:AD-LaGa}
  Given dual profunctors
  \[ f : P \promap Q, \quad g = Df : Q \promap P. \]
Let  $(I,F)$ be the cut for  $Q \times P^\op$ associated to $f$, and
$(J,G)$ the cut for $P \times Q^\op$ associated to
$g$.
\begin{itemize}
\item[a.] $(J,G) = (F^\op,I^\op)$.
\item[b.] $\Gamma G = (\Lambda I)^\op, \quad \Lambda J = (\Gamma F)^\op$.
\item[c.] $\Gamma g = (\Lambda f)^\op, \quad \Lambda g = (\Gamma f)^\op$.

\end{itemize}
\end{corollary}

\begin{proof}
Part a is by Lemma \ref{lem:posPQ-Dfq}.   
Parts b and c are then immediate from a.
\end{proof}

\subsection{Extending $\Lambda$ and $\Gamma$
to the next level}
We have looked at cuts $(I,F)$ for $Q \times P^\op $.
Proceeding to the next level, we look at cuts
$(\cI,\cF)$ for $\Pro(P, Q) = (Q \times P^\op)^\dual$.
Elements in this latter set are cuts $(I,F)$ for $Q \times P^\op$
partially ordered by $(I,F) \leq (I^\prime, F^\prime)$ if $I \sus I^\prime $.
Thus the down-set $\cI $ consists of cuts $(I,F)$ closed under
taking cuts with smaller $I$'s. Similarly the up-set $\cF$ is
closed under taking larger $I$'s (and so smaller $F = I^c$'s).

\medskip
We have a map
\begin{equation} \label{eq:AD-Lambda}
  \Lambda : U\Pro(P, Q)  \pil  (UQ \times UP^\op)^\dual, \quad
  f \mapsto  (\Lambda f,-)
\end{equation}
and a map
\begin{equation} \label{eq:AD-Gamma}
\Gamma : U\Pro(P, Q)  \pil  (UQ \times UP^\op)^\dual, \quad
  f  \mapsto  (-,\Gamma f).
\end{equation}
By \eqref{eq:PQ-gu} the map $\Lambda$ induces an isotone map of posets:
\begin{equation*}
  \Lambda_U^\ii : \Pro(P,Q)^\dual \pil  (UQ \times UP^\op)^{\doublehat{}},
  \quad (\cI,\cF)  \mapsto  (-, \Lambda(\cF)^\uparrow).
\end{equation*}
So the image of the cut
$(\cI,\cF)$ is the cut in the Boolean lattice
$(UQ \times UP^\op)^\dual$ whose up-set
is $\Lambda(\cF)^\uparrow$,
the up-set generated by all $\Lambda f$ for $f \in \cF$.

The map $\Gamma $ induces an isotone  map of posets:
\begin{equation*}
  \Gamma_U^!: \Pro(P, Q)^{\dual}  \pil
  (UQ \times UP^\op)^{\doublehat{}},
  \quad (\cI, \cF)  \mapsto  (\Gamma(\cI)^{\downarrow},-).  
\end{equation*}
So the image of the cut $(\cI,\cF)$ is the cut in the Boolean lattice
$(UQ \times UP^\op)^\dual$ whose ideal is
$\Gamma(\cI)^{\downarrow}$, the ideal generated by
the {\em complements} (see \eqref{eq:AD-Gamma})
of all $\Gamma f$ for $f \in \cI$,

\subsection{The main theorem: Preserving the cut}

In order for the left and right boundaries of a cut $(I,F)$ to give
enough information we need to make sure that minimal elements
of up-sets of $P$ and $Q$ exists. 
A poset $P$ is {\it well-founded} if every subset of $P$ has a minimal element.
Equivalently, any descending chain of elements
  \[ p_1 \geq p_2 \geq \cdots \geq p_n \geq \cdots
  \]
  stabilizes, i.e. for some $N$ we have $p_n = p_N$ for $n \geq N$.

The following theorem is a strong generalization of the results 
in several articles \cite{EHM}, \cite{FGH}, \cite{F-LP},\cite{JuKaMa}, see
Section \ref{sec:AD-LP} for more on this.
The most significant tool in the argument is Zorn's lemma (which is
equivalent to the axiom of choice). Note also that $\Pro(P,Q)$
is a complete distributive lattice and so has all joins (colimits)
and meets (limits).
\begin{theorem}[Preserving the cut] \label{thm:AD-duality}
  Let $P$ and $Q$ be well-founded posets,
  and $(\cI,\cF)$ a cut for
  $\Pro(P, Q)$. Then $(\Gamma(\cI)^\downarrow,\Lambda(\cF)^\uparrow)$ is a cut
  for the Boolean lattice $(UQ \times UP^\op)^\dual$. In other
  words, the maps $\Gamma_U^! = \Lambda_U^\ii$.
\end{theorem}

\begin{example} \label{eks:AD:53}
  Consider profunctors $\Hom([5],\di{[3]})$ where
  $\di{[3]} = [3] \cup \{\infty\}$. Let $\cI$ be the ideal
  consisting of all $f$ not taking the value $\infty$. The graph of such
  an $f$ is shown in Figure \ref{fig:gr-as}. Such a path
  in the rectangle $U[3] \times U[5]^\op$ is a {\em right path}.
  Let $\cF$ be the complement
  up-set, consisting of all $g$ which take the value $\infty$ for some
  argument in $[5]$. The ascent of such a $g$ is also shown in Figure
  \ref{fig:gr-as}. Such a path
  in the rectangle $U[3] \times U[5]^\op$ is an {\em up path}.
  Theorem \ref{thm:AD-duality} above says that given any subset $S$ of
  $[5] \times [3]$, exactly one of the following holds:
  i. $S$ contains an up path, ii. the complement $S^c$ contains a right path.
  An earlier observation of this is in \cite[Lemma 4.4] {FV-Bi}.
\end{example}

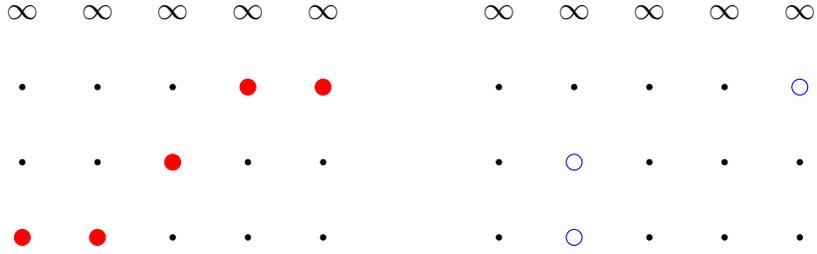
\begin{figure}
\begin{tikzpicture}
  {    
\filldraw[color=red] (1,1) circle (3pt);
\filldraw[color=red] (2,1) circle (3pt);
\filldraw[color=black] (3,1) circle (1pt);
\filldraw[color=black] (4,1) circle (1pt);
\filldraw[color=black] (5,1) circle (1pt);

\filldraw[color=black] (1,2) circle (1pt);
\filldraw[color=black] (2,2) circle (1pt);
\filldraw[color=red] (3,2) circle (3pt);
\filldraw[color=black] (4,2) circle (1pt);
\filldraw[color=black] (5,2) circle (1pt);

\filldraw[color=black] (1,3) circle (1pt);
\filldraw[color=black] (2,3) circle (1pt);
\filldraw[color=black] (3,3) circle (1pt);
\filldraw[color=red] (4,3) circle (3pt);
\filldraw[color=red] (5,3) circle (3pt);

\draw node at (1,4) {$\infty$};
\draw node at (2,4) {$\infty$};
\draw node at (3,4) {$\infty$};
\draw node at (4,4) {$\infty$};
\draw node at (5,4) {$\infty$};
}
\end{tikzpicture}
\hskip 1.5cm
\begin{tikzpicture}
  {
\filldraw[color=black] (1,1) circle (1pt);
\draw[color=blue] (2,1) circle (3pt);
\filldraw[color=black] (3,1) circle (1pt);
\filldraw[color=black] (4,1) circle (1pt);
\filldraw[color=black] (5,1) circle (1pt);

\filldraw[color=black] (1,2) circle (1pt);
\draw[color=blue] (2,2) circle (3pt);
\filldraw[color=black] (3,2) circle (1pt);
\filldraw[color=black] (4,2) circle (1pt);
\filldraw[color=black] (5,2) circle (1pt);

\filldraw[color=black] (1,3) circle (1pt);
\filldraw[color=black] (2,3) circle (1pt);
\filldraw[color=black] (3,3) circle (1pt);
\filldraw[color=black] (4,3) circle (1pt);
\draw[color=blue] (5,3) circle (3pt);

\draw node at (1,4) {$\infty$};
\draw node at (2,4) {$\infty$};
\draw node at (3,4) {$\infty$};
\draw node at (4,4) {$\infty$};
\draw node at (5,4) {$\infty$};
}
\end{tikzpicture}
\caption{Profunctors from $[5]$ to $[3]$. A graph to the left and an
 (unrelated) ascent to the right.}
\label{fig:gr-as}
\end{figure}

Before proving Theorem \ref{thm:AD-duality} we state two alternative
formulations of this theorem.

\begin{theorem} \label{thm:AD-ad}
  Let $P$ and $Q$ be well-founded posets,
  and $(\cI,\cF)$ a cut for
  $\Pro(P, Q)$. The
  up-set  $\Lambda^\ii(\cF)$ for $(UQ \times UP^\op)^\dual$
  generated by all $(\Lambda f,-)$ for $f \in \cF$, and
  the up-set $(\Gamma^\op)^\ii(\cI^\op)$ for $(UP \times UQ^\op)^\dual$
  generated by all $((\Gamma f)^\op,-)$ for $f \in \cI$, are
  Alexander dual up-sets.
  (This is the version applied in Stanley-Reisner theory, giving
  Alexander dual monomial ideals, see Section \ref{sec:AD-LP}.)
\end{theorem}

\begin{proof}
By the above Theorem \ref{thm:AD-duality}, 
$(\Gamma(\cI)^{\downarrow})^\op = \Gamma^!(\cI)^\op$ is the  Alexander dual
up-set of the up-set $\Lambda^\ii(\cF) = \Lambda(\cF)^\uparrow$.
But $\Gamma^!(\cI)^\op = \Gamma^{! \, \op}(\cI^\op)$ and
by \eqref{eq:gud} $(\Gamma^{\op})^\ii$ equals $(\Gamma^{!})^\op$. 
  \end{proof}

By Corollary \ref{cor:AD-LaGa} we have a commutative diagram
\begin{equation} \label{eq:AD-LGdiagram} \xymatrix{
     U\Pro(P,{Q}) \ar^{\Lambda}[r] \ar^{D}[d]
      & (UQ \times UP^\op)^{\dual} \ar^{D}[d]\\
      U\Pro(Q,{P}) \ar^{\Gamma}[r] & (UP \times UQ^\op)^{\dual}
    }.
\end{equation}
In the proof below we write $\Lambda_{P,Q}$ and $\Gamma_{Q,P}$ for these
horizontal maps.

  \begin{theorem} \label{thm:AD-lala}
    Let $P$ and $Q$ be well-founded posets. The following diagram
    commutes:
     \[ \xymatrix{
      \Pro(P,{Q})^\dual \ar^{\Lambda_U^\ii}[r] \ar^{D}[d]
      & (UQ \times UP^\op)^{\doublehat{}} \ar^{D}[d]\\
      \Pro(Q,{P})^\dual \ar^{\Lambda_U^\ii}[r] & (UP \times UQ^\op)^{\doublehat{}}
    }.
  \]
\end{theorem}

\begin{proof}
  By the commutative diagram \eqref{eq:AD-LGdiagram},
  $\Gamma_{Q,P}^\op = \Lambda_{P,Q}$. Switching $P$ and $Q$ we have
  $\Gamma_{P,Q}^\op = \Lambda_{Q,P}$. Hence $(\Gamma_{P,Q}^\op)_U^\ii$ equals
  $(\Lambda_{Q,P})_U^\ii$. By the above Theorem \ref{thm:AD-ad}
$(\Gamma_{P,Q}^\op)_U^\ii$ is the dual of $(\Lambda_{P,Q})_U^\ii$ and so 
  $(\Lambda_{Q,P})_U^\ii$ is the dual of $(\Lambda_{P,Q})_U^\ii$.
\end{proof}

\begin{problem} \label{pro:AD:functorial}
  Theorem \ref{thm:AD-duality}
  may be seen as a map
  \[ \Pro(P,Q)^\dual \pil \Pro(UP, UQ)^\dual. \]
  However considering the natural maps $UP \pil P$ and
  $UQ \pil Q$ none of the natural functorial ways to get a map
  as above, gives the map $\Gamma^! = \Lambda^{\ii}$. For instance
  we can not see a natural factorization of the above map
  through $\Pro(UP, Q)^\dual$. 
  
  Is there a way to understand or formulate the theorem  to get
  a functorial construction, which works for maps $P^\prime \pil P$
  and $Q^\prime \pil Q$ (or in some other way)?
\end{problem}

\begin{proof}[Proof of Theorem \ref{thm:AD-duality}.]
  We need to show that the ideal $\Gamma(\cI)^{\downarrow}$ is the complement
  of the up-set $\Lambda(\cF)^{\uparrow}$ in the Boolean lattice $(UQ \times UP^\op)^\dual$.

  \medskip
\noindent {\bf Part I.} We first show that
  for any $(I,F) \in  \cI $ and $(J,G) \in \cF$ that 
  \[ \Gamma F \cap \Lambda J  \neq \emptyset . \]
  By Lemma \ref{lem:posP-IFd}
  this shows that $\Gamma(\cI)^{\downarrow}$ and $\Lambda(\cF)^{\uparrow}$
  are disjoint.
  
  If $(I,F) \in \cI$ and $(J,G) \in \cF$ then $(F \cap J) \sus
  Q \times P^\op$ is $\neq \emptyset$ by
  Lemma \ref{lem:posP-IFd}. Let $(q_0,p_0) \in F \cap J$.
  Define $(q_i,p_i)$ in $F\cap J$
  successively as follows. If there is a $q < q_i$
  such that $(q,p_i)$ is in $F$, then let $(q_{i+1},p_{i+1})$  
  be $(q,p_i)$ which is automatically also in $J$.
  If there is $p < p_i$ such that $(q_i,p)$ is also in $J$,
  then let $(q_{i+1},p_{i+1})$ be $(q_i,p)$ which is automatically also
  in $F$. We get chains $p_0 \geq p_1 \geq \cdots $ and
  $q_0 \geq q_1 \geq \cdots$. These must stabilize and for the stable
  value pair $(q,p)$ we cannot continue the construction and then this
  pair will be in $\Gamma F \cap \Lambda J$. 

\medskip
\noindent {\bf Part II.}
Let $S \sus UQ \times UP^\op$ be such that
$S \cap \Lambda J \neq \emptyset$ for {\it every} $(J,G) \in \cF$.
This means that $(S^c,S)$ is not in $\Lambda(\cF)^{\uparrow}$. 
We show that $S \supseteq \Gamma F$ for {\it some} $(I,F) \in \cI$
and so $(S^c,S)$ is in $\Gamma(\cI)^{\downarrow}$. This shows that the
union of $\Gamma(\cI)^{\downarrow}$ and $\Lambda(\cF)^{\uparrow}$ is all
of $UQ \times UP^\op$.

Let $T$ consist of $f \in \Pro(P,Q)$ such that $\Lambda f \cap S =
\emptyset$. Clearly the minimum $\bfnu$ of $\Pro(P,Q)$ is in $T$, so
$T$ is non-empty. We use Zorn's lemma to show that $T$ has a maximal element.
So let
\begin{equation} \label{eq:AD-fkjede}
  f_1 \leq f_2 \leq \cdots \leq f_n \leq \cdots
  \end{equation}
  be a chain of functions in $T$, and $f$ be the colimit (join) of these.
  We claim that $\Lambda f \cap S = \emptyset$. Let
$(q,p) \in \Lambda f$, so $ q \in f(p) $ but $q \not \in f(p^\prime)$ for
every $p^\prime < p$. Then $q \in f_i(p)$ for some $i$ and
$q \not \in f_i(p^\prime)$ since $f_i(p^\prime) \sus f(p^\prime)$. Then
$(q,p) \in \Lambda f_i$ and so is not in $S$. Thus $f$ is an upper bound
in $T$ for the chain in \eqref{eq:AD-fkjede}.
By Zorn's lemma $T$ has a maximal element $f$, and take note that
$f$ must be in $\cI$, by the requirement on $S$.

\begin{claim}
  $\Gamma f \sus S$.
\end{claim}

  By the second and third lines of Part II, this
  finishes the proof.

\medskip
  
\noindent {\it Proof.} Suppose not, so there is
  $(q_0,p_0) \in \Gamma f \backslash S$.
Recall $q_0$ is minimal in the complement $(f(p_0))^c$. Let
\[ \tilde{f} (p) = \begin{cases} f(p) \cup \{ q_0 \}, &
    p \geq p_0, \\ f(p) & \text{otherwise}.
  \end{cases} \]
This is a profunctor $\tilde{f} : P \pil Q$. Let us show
\[ \Lambda \tilde{f} \sus \Lambda f \cup\{ (q_0,p_0)\}.\]
That $(q,p) \in \Lambda \tilde{f}$ means that
  $q \in \tilde{f}(p)$ and
  $q \not \in \tilde{f}(p^\prime)$ for $p^\prime < p$.
  \begin{itemize}
\item[1.]  If not $p \geq p_0$ then $\tilde{f}(p^\prime)  =
    f(p^\prime)$ for all
  values $p^\prime \leq p$. Therefore $(q,p)$ is in $\Lambda f$.
\item[2.]  Let $p \geq p_0$. Then either i. $q \in f(p)$ or ii. $q = q_0$.
  In case i. $(q,p) \in \Lambda f$. In case ii, where
  $q = q_0$, 
  if $p > p_0$ then $(q_0,p)$ would  not be in $\Lambda \tilde{f}$ since
  $q \in \tilde{f}(p)$ and $q \in \tilde{f}(p_0)$. 
  We are thus left with $(q,p) = (q_0,p_0)$.
\end{itemize}
As a consequence we also have $\Lambda \tilde{f} \cap S = \emptyset$.
But by definition of $S$ then $\tilde{f} \in T$. This contradicts
$f$ being maximal in $T$. Therefore no such $(q_0,p_0)$ could exist,
and $\Gamma f \sus S$.

\end{proof}

\section{Applications to finite posets and Stanley-Reisner
  ideals} \label{sec:AD-LP}

In \cite{EHM} they show that letterplace ideals $L([n],P)$ and
co-letterplace ideals $L(P,[n])$ are Alexander dual ideals, and this was
taken considerably further in \cite{FGH}.

Here we define such square free ideals in the full general setting
for finite posets $P$ and $Q$ and cuts in $\Pro(P,Q)$. We show
how the result of \cite{EHM} above is a consequence.
In the following for ease of notation we denote a polynomial
ring $k[x_{UQ \times UP^{op}}]$ as $k[x_{Q \times P^\op}]$. 

\begin{definition} Let $P$ and $Q$ be {\em finite} posets.
  From the cut $(\cI,\cF)$ for $\Pro(P,{Q})$,
  we get the cut $(\Gamma(\cI)^\downarrow, \Lambda(\cF)^\uparrow)$ for
  $(UQ \times UP^\op)^\dual$.

  \medskip
  \noindent{\em The $\Lambda$-ideal.}
  Since $UQ \times UP^\op$ is simply a set, this gives a
  Stanley-Reisner ideal in $k[x_{Q\times P^\op}]$.
  We denote it $L_\Lambda(\cF;P,Q)$, or simply $L_\Lambda(\cF)$. As $f$ varies in $\cF$,
  it is generated by the
  $\Lambda f$ (or rather the squarefree monomials
  $ \prod_{(q,p)\in \Lambda f}x_{q,p}$).

\medskip
  \noindent{\em The $\Gamma$-ideal.}
  The dual cut $(\Lambda(\cF)^{\uparrow \op}, \Gamma(\cI)^{\downarrow \op})$ for
  $(UP \times UQ^\op)^\dual$, also corresponds to 
  a Stanley-Reisner ideal in $k[x_{P \times Q^\op}]$. We denote this ideal as
  $L_\Gamma(\cI;P,Q)$, or simply $L_\Gamma(\cI)$. As $f$ varies in $\cI$,
  it is generated by the $(\Gamma f)^{\op}$ (or rather the squarefree monomials
  $\prod_{(p,q) \in (\Gamma f)^{\op}} x_{p,q}$).
\end{definition}

By Theorem \ref{thm:AD-ad} above, the $\Lambda$-ideal and $\Gamma$-ideal
are Alexander dual ideals.
  

\begin{remark} Let $Q = [n]$, the chain on $n$ elements,
and $\cI \sus \Pro^{< \infty}(P,{[n]})$.  The $\Lambda$-ideals
  $L_\Lambda(\cF;P,[n])$ are the letterplace ideals of \cite{FGH} and
  are shown to be Cohen-Macaulay ideals. The Alexander dual $\Gamma$-ideals
  $L_\Gamma(\cI;P,[n])$ are the co-letterplace ideals in loc.cit. and
  thus have linear resolutions, \cite[Thm.5.56]{Mi-St}.
  The above definition and Theorem \ref{thm:AD-ad}
  may thus be seen as a full generalization of the setting of \cite{FGH}.
\end{remark}

\begin{remark} With the same setting as in the above remark, 
  the $\Lambda$-ideals
    $L_\Lambda(\cF;P,[n])$ define simplicial balls by the Stanley-Reisner
    correspondence, \cite{DFN}. Furthermore there is a very simple
    description of the Stanley-Reisner ideal of their boundaries, which
    are simplicial spheres. This gives the construction of an enormous
    amount of simplicial spheres, due to the freedom in choosing
$P,[n]$ and $\cI$.
\end{remark}

Corollary  \ref{cor:PQ-ad} has the following consequence, see also Example
\ref{eks:AD:53}.

  \begin{corollary} \label{cor:AD-LP:infty} If $P$ or $Q$ has a maximal element,
    the ideals $L_\Gamma(\Pro^{< \infty}(P, {Q}))$ and
    $L_\Gamma(\Pro^{< \infty}(Q,{P}))$ are Alexander dual ideals.
  \end{corollary}

    \begin{definition}
      Given an isotone map $f : P \pil Q$, let the graph
      \[ \Gamma_i f = \{ (f(p),p^\op)\, | \, p \in P \} \sus UQ \times UP^\op. \]
      This gives a map
      \[ \Gamma_i : \Hom(P,Q) \pil (UQ \times UP^\op)^\dual, \quad
        f \mapsto (-, \Gamma_i f), \]
      and so
      \[ \Gamma_i^\op : \Hom(P,Q)^\op \pil (UP \times UQ^\op)^\dual, \quad
         f^\op \mapsto ((\Gamma_i f)^\op, -), \]
      The {\em isotonian} ideal $L(P,Q)$ is the ideal in
      $k[x_{P \times Q^\op}]$ generated by the monomials
      $\prod_{(p,q) \in (\Gamma_i f)^\op} x_{(p,q)}$ as $f$ varies in $\Hom(P,Q)$.
      If $P$ is the chain $[n]$, it is called a letterplace ideal.
      If $Q = [n]$, it is a co-letterplace ideal, \cite{FGH}.
    \end{definition}

    The canonical map $Q \pil \di{Q}$ induces an isotone map
    \begin{equation} \label{eq:AD-alpha}
      \alpha : \Hom(P,Q) \pil \Hom(P,\di{Q}) = \Pro(P,Q), \quad
      g  \overset{\alpha}{\mapsto} (p {\, \mapsto} \,  \widehat{g(p)} =
      (-, \upa g(p))),
      \end{equation}
      and by Subsection \ref{subsec:PQ-adj}
      \begin{equation} \label{eq:AD-alphaP}
      \alpha^{!} : \Hom(P,Q)^\dual \pil \Pro(P,Q)^{\dual}.
    \end{equation}
    We get commutative diagrams:
    \begin{equation} \label{eq:ad-hompq}
      \xymatrix{ U\Hom(P,Q)^{\op} \ar[r]^{\alpha^\op}
\ar[rd]_{\Gamma_i^\op}
    & 
    U\Pro(P,Q)^{\op}\, ,  \ar[d]^{\Gamma^{\op}} \\
    & (UP \times UQ^\op)^\dual}
  \quad
  \xymatrix{ \Hom(P,Q \widehat{)^{\op}} \ar[r]^{(\alpha^\op)^\ii}
\ar[rd]_{(\Gamma_{i,U}^\op)^{\ii}}
    & 
    \Pro(P,{Q} \widehat{)^{\op}}  \ar[d]^{(\Gamma_U^\op)^\ii = (\Gamma_U^!)^\op} \\
    & (UP \times UQ^\op)^{\doublehat{}  }}
  \end{equation}
  Note that the isotonian ideal $L(P,Q)$ is the squarefree monomial ideal
  associated to the image by $(\Gamma_i^\op)^\ii$ of the cut
  $(\emptyset, \Hom(P,Q)^\op)$ for $\Hom(P,Q \widehat{)^{\op}}$.

The poset  $P$ is a {\it forest} if for every two incomparable $p_1$ and $p_2$
    in $P$, there is no $p \in P$ such that $p \leq p_1$ and $p \leq p_2$. The
    connected components of the Hasse diagram of $P$ are then trees with
    the roots on top.

    
    \begin{lemma}
      If $P$ is a forest, then:
      \begin{itemize} 
      \item[a.] For the map in \eqref{eq:AD-alphaP},
        $\alpha^{!}(\Hom(P,Q))$ is the down-set
        $\Pro^{< \infty}(P,Q)$.
      \item[b.] In the diagram \eqref{eq:ad-hompq}
        \[ (\Gamma_{i,U}^{\op})^\ii ((\emptyset, \Hom(P,Q)^\op)) =
          (\Gamma_U^{\op})^\ii((-, \Pro^{< \infty}(P,Q)^\op)). \]
        \end{itemize}
Hence $L(P,Q) = L_\Gamma(\Pro^{< \infty}(P,Q))$. 
    \end{lemma}

    \begin{proof}
By the right diagram of \eqref{eq:ad-hompq}, part a above implies part b.
So we do  part a. 
      We show that if $f$
      is a profunctor in $\Pro^{< \infty}(P,Q)$,
      then there exists a $g : P \pil Q$ such that
      $\Gamma f \supseteq \Gamma_i g$.
      This implies that $\alpha(g) \geq f$ giving part a. It also implies
      part b.
      \begin{itemize}
      \item[1.] If $p$ is maximal in $P$ then let $q$ be a minimal element
      in the non-empty set $f(p)^c$. Define then $g(p) = q$.
\item[2.] Suppose $p$ is not maximal. Suppose $g(p^\prime)$ is defined for
      $p^\prime > p$, such that
      $(g(p^\prime), p^{\prime \, \op}) \in \Gamma f$.
      Let $p^\prime $ be the unique cover
      of $p$ (since $P$ is a forest). Then $f(p) \sus f(p^\prime)$ and so $f(p)^c \supseteq
      f(p^\prime)^c$. Given $g(p^\prime) \in f(p^\prime)^c$. Then there will exist
      a minimal $q$ in $f(p)^c$ with $q \leq g(p^\prime)$. Then define
      $g(p) = q$. This gives the map $g$.
    \end{itemize}
    \end{proof}

We specialize to $Q = [n]$ and have the consequence:

    \begin{corollary}
      The letterplace ideal $L([n],P)$ and the co-letterplace ideal
      $L(P,[n])$ are Alexander dual ideals.
      \end{corollary}
 
    \begin{proof}
      The ideal $L(P,[n])$ is  $L_\Gamma(\Pro^{< \infty}(P,[n]))$ since
      $\di{[n]} = [n] \cup \{ \infty \}$. Its
      Alexander dual is the ideal $L_\Gamma(\Pro^{< \infty}([n],P)$
      by Corollary \ref{cor:AD-LP:infty}.
      By the above lemma this identifies as $L([n],P)$.
    \end{proof}
    In \cite{HeQuSh} they fully characterize for which $P$ and $Q$
    the isotonian ideals $L(P,Q)$ and
    $L(Q,P)$ are Alexander dual.

\section{Topology on $\Pro(P,Q)$} \label{sec:top}

We want to get a setting where $P$ and $Q$ may be infinite but
the ascents $\Lambda f$ and graphs $\Gamma f$ are finite.
Then we get $\Lambda$- and $\Gamma$-ideals in infinite-dimensional
polynomial rings.
This can be achieved if we
put a suitable topology on $\Hom(P,Q)$. 

\subsection{Defining the topology}

We define a topology on $\Pro(P,Q)$ by defining a basis of
open subsets.

\begin{definition} \label{defi:PosU} Given a large profunctor $\ov{f}$ in
  $\Pro^L(P,Q)$ and a small profunctor $\ul{f}$ in $\Pro_S(P,Q)$.
  The set of profunctors $f$ such that $\ul{f} \leq f \leq \ov{f}$
  is denoted $U(\ul{f},\ov{f})$.
\end{definition}

The following is immediate.

\begin{lemma}
  Given two such sets $U_1 = U(\ul{f}_1, \ov{f}_1)$ and
  $U_2 = U(\ul{f}_2, \ov{f}_2)$, let $\ul{f}$ and $\ov{f}$ be the profunctors
  defined by the join and meet
  \[  \ul{f} = \ul{f}_1 \vee \ul{f}_2, \quad
    \ov{f} = \ov{f}_1 \wedge \ov{f}_2. \]
  Then
  \[ U_1 \cap U_2 = U(\ul{f}, \ov{f}). \]
\end{lemma}


As a consequence the sets $U(\ul{f}, \ov{f})$ form the basis of a topology
on $\Pro(P,Q)$. Note that if $P$ and $Q$ are finite we get the discrete
topology on $\Pro(P,Q)$.

\begin{observation} \hskip 2mm
  a. There is a one-one correspondence between
  open down-sets in $\Pro(P,Q)$ and down-sets in
  $\Pro^L(P,Q)$. 

  b. There is a one-one correspondence between
  open up-sets in $\Pro(P,Q)$ and up-sets in
  $\Pro_S(P,Q)$. 
\end{observation}

Note also by Lemma \ref{lem:settD} that the open set $U(\ul{f}, \ov{f})$
is mapped by the dual $D$ to the open set $U(D\ov{f}, D\ul{f})$. Hence
we get:

\begin{lemma} The duality map $D$ of Lemma \ref{lem:settD} is a homeomorphism
of topological spaces.
\end{lemma}



\subsection{Generalities on topology}

In any topological space $X$ there is distinguished type of open subsets:
those that are the interiors of closed subsets, called {\it regular
  open sets}. In fact, let
$ \op X$ be the poset of open subsets of $X$, and $\cl X$ the poset
of closed subsets of $X$. There is a Galois correspondence:
\begin{center}
\begin{eqnarray*}
 &\op X  \bihom{\overline{()}}{()^\circ}  \cl X,  &\\
& U  \mapsto \ov{U} = \text{ closure of } U, \quad
C   \mapsto C^{\circ} =  \text{interior of } C& 
\end{eqnarray*}
\end{center}

The fix points in $\op X$ for the Galois connection, i.e. the fix points
in $\op X$ for the composition of closure $\ov{{}}$ and interior ${}^\circ$,
are the regular open sets. We denote these as $\rego X$, the
open subsets which are interiors of closed subsets.

\begin{lemma}
  The map on open subsets $U \mapsto \ov{U}^c$ where $c$ denotes complement,
  gives an involution $ \rego X \mto{i} \rego X$. Furthermore
  $\ov{U}^c = (U^c)^\circ$. 
\end{lemma}

\begin{proof}
  Let $V = \ov{U}^c$. We will show that $\ov{V}^c = U$ by proving
  inclusion either way.

  1. We get $V \sus U^c$ and since the latter is closed, $\ov{V} \sus U^c$
  and so $\ov{V}^c \supseteq U$.
  
  2. We also have $\ov{U}^c \sus \ov{V}$ and so $\ov{U} \supseteq \ov{V}^c$.
  Since $U$ is the interior of $\ov{U}$ we get $U \supseteq \ov{V}^c$.
  
  3. We now show that $V$ is regular. Let $W$ be an open subset of $\ov{V}$.
  Then $W^c \supseteq \ov{V}^c = U$, so $W^c \supseteq \ov{U} = V^c$.
  Thus $W \sus V$. Hence $V$ is the largest open subset of $\ov{V}$, and
  so is the interior.

  4. Since $U = \ov{V}^c$ we get $U^c = \ov{V}$, and so
  \[ (U^c)^\circ = \ov{V}^\circ = V = \ov{U}^c. \]
\end{proof}

\subsection{Open and closed down-sets}

\begin{lemma} \label{lem:appCO} Let $\cI$ be a down-set in $\Pro(P,Q)$.
  The closure $\ov{\cI}$ and the interior $(\cI)^\circ$ are also down-sets.
  By duality the analog also holds for up-sets.
\end{lemma}

\begin{proof}
  Let $F \in \ov{\cI}$ and let $G \leq F$. We will show $G$ is also in
  $\ov{\cI}$. If not there is an open $U(b,c)$ disjoint from $\cI$
  with $G \in U(b,c)$, so $b \leq G \leq c$.

  Then $F \geq b$ with $b$ small and not in $\cI$. Then
 $F \in U(b,\bfen)$ with this basis open set disjoint from $\cI$
  since $b$ is not in $\cI$ and $\cI$ is a down-set. But then $F$ is not
  in the closure of $\cI$, contrary to assumption.

  As for the interior, if $\cF$ is the complement up-set of $\cI$, then
  by a similar argument $\ov{\cF}$ is a up-set, and so
  the complement $(\ov{\cF})^c = \cI^\circ$ is a down-set.
\end{proof}

\begin{lemma} \label{lem:top-gap}
Let $\cI$ be an open down-set $\cI$. Then 
$\cI$ is regular if and only if $\overline{\cI}\backslash \cI$ does not contain large
maps.

Similarly an open up-set $\cF$ is regular if and only if $\overline{\cF}\backslash
\cF$ does not contain small maps.

\end{lemma}

\begin{proof} We only show the first statement. Suppose $\cI$ is regular. 
 Let $f$ in $\ov{\cI}\backslash \cI$ be large.
 The closure $\ov{\cI}$ is a down-set.
 It contains $f$ and so the open set
  $U(\bfnu,f)$, which must then be in the interior of $\ov{\cI}$ which is
  $\cI$. This contradicts $f$ being in the gap. Similarly one can argue
  that an $f$ in the gap cannot be bounded, by considering the up-set $\cF$.

  If $\cI$ is not regular, there is some open
  $U = U(\underline{f}, \overline{f})
  \sus \overline{\cI}$ with $U \not \sus \cI$. But then
  $\overline{f} \in \overline{\cI} \backslash \cI$, and so
  $\overline{\cI} \backslash \cI$  contains $\overline{f}$ which is large.
\end{proof}


\begin{definition} A {\em Dedekind cut} in the poset $\Pro(P,Q)$ 
is a pair $[\cI,\cF]$ where $\cI$ is an regular open down-set of
    $\Pro(P,Q)$ and 
    $\cF$ is its image by the involution $i$
    (so $\cF$ is an regular open up-set).
The {\em gap} of the Dedekind cut is
\[ \cG = \Pro(P, Q) \backslash (\cI \cup \cF).  \]
\end{definition}

\begin{corollary} Let $\cI$ be an open down-set and $\cF$
  a disjoint open up-set. Then $\cI$ and $\cF$ form a Dedekind
  cut if and only if the gap 
  $\Pro(P, Q) \backslash (\cI \cup \cF)$ between them, does not
  contain large or small maps.
\end{corollary}

\begin{remark}\label{rem:top-real}
  For the rational numbers $\QQ$ the distributive lattice $\di{\QQ}$ consists
  of $\pm \infty$, the irrational numbers, and the rational numbers
  doubled, since a rational number $q$ gives two cuts:
  \[ (\langle -\infty,q),[q,+\infty \rangle), \quad
    (\langle -\infty,q],(q,+\infty \rangle). \]
  However if we have the natural topology on $\QQ$, Dedekind cuts
  give the real numbers together with $\pm \infty$. 
  \end{remark}

\section{Profunctors between natural posets}
\label{sec:natural}

We introduce the class of natural posets as a suitable generalization
of the poset of natural numbers. We then have good criteria for when open
down-sets of $\Pro(P,Q)$ are closed or regular.
We show, Theorem \ref{thm:natural:fg}, that an open down-set
is clopen (closed and open) if and only if it has a finite number of maximal elements.
The purpose of this section is to get
in Section \ref{sec:monwell}
a version of Theorem
\ref{thm:AD-duality} which applies to construct Alexander dual
ideals in infinite-dimensional polynomial rings.
Our main application is to profunctors
$\NN \promap \NN$, see \cite{Fl-St}, and briefly indicated in the
last Section \ref{sec:mon}.


\subsection{Well partially ordered sets}
A poset $P$ is {\it well partially ordered} if the following two conditions
holds:

\begin{itemize}
\item[i.] Any descending chain in $P$ stabilizes.
\item[ii.] Any antichain in $P$ is finite.
\end{itemize}

The following characterization of well partially ordered sets is classical,
see for instance \cite{Kr}.

\begin{proposition} The following are equivalent for a poset $P$.

  \begin{enumerate}
  \item $P$ is well partially ordered.
  \item Any up-set of $P$ is finitely generated.
  \item The ascending chain condition holds for up-sets in $P$,
    equivalently the descending chain condition holds for $\di{P}$.
    \end{enumerate}
  \end{proposition}

\begin{lemma} \label{lem:AD-fin}
Let $P$ and $Q$ be well partially ordered sets, and $f : P \promap {Q}$ an
isotone map.
\begin{itemize}
\item[a.] If $f$ is large, then $\Gamma f$ is finite.
\item[b.] If $f$ is small, then $\Lambda f$ is finite.
  \end{itemize}
\end{lemma}

\begin{proof}
a. Let $f$ have profile $(I,F)$. The projection of $\Gamma f$ onto
$P$ will be contained in $I$ since the complement $f(p)^c$ is empty
for $p \in F$. Furthermore for $p \in I$ the up-set $f(p)^c$
has only a finite set of minimal elements by the above proposition.

b. Let $f$ have co-profile $(J,G)$. The projection of $\Lambda f$ onto
$Q$ is contained in $J$. If $q \in J$, the set of $p \in P$ such that
$q \in f(p)$ is a up-set in $P$ and so has a finite set of minimal elements.
Hence $\Lambda f$ is finite.
\end{proof}

\subsection{Natural posets}
The class of well partially ordered sets is very large, it includes all well-ordered
sets. For our purposes we need an extra condition so that our
posets are more like natural numbers.

\begin{definition} A poset $P$ is {\em down finite} if for each
  $p$ in $P$, the principal down-set $\doa p$ is finite.
  A {\em \napo} poset (in analogy with natural numbers) is a poset which is
  well partially ordered and down finite.
\end{definition}


\begin{lemma}
  Suppose every antichain in $P$ is finite. Then $P$ is \napo iff
  there exists a chain of finite down-sets in $P$
  \[ I^0 \sus I^1 \sus \cdots \sus I^n \sus \cdots \]
  such that $\cup_{i \geq n} I^n$ equals $P$.
\end{lemma}

\begin{proof} The if statement is clear. Suppose then $P$ is natural.
  For $p \in P$ the principal ideal $\doa p$ is finite. Let the {\em height} of
  $p$, denoted $\ell(p)$, be
  the length of the longest chain in $\doa p$. Then let $I^r$ be the set of
  elements of $P$ of height $\leq r$. Then $I^r\backslash I^{r-1}$ is the
  set of elements of height exactly $r$. These form an anti-chain in $P$
  and so a finite set. Hence each $I^r$ is finite. Clearly the union of
  the $I^r$ is all of $P$.
\end{proof}

\begin{lemma} \label{lem:top-single} Let $Q$ be down finite.
Then for any large $f$ in $\Pro(P,Q)$ the
  open set $U(\bfnu,f)$ is also closed.
\end{lemma}

\begin{proof}
  Let $f$ have profile $(I,F)$.
  Consider the set $X$ of all pairs $(p,x)$ where $p \in I$ and $x$
  is minimal in $f(p)^c$. For each pair $(p,x) \in X$, let $g_{px}$ be
  the smallest isotone map $g : P \pil \di{Q}$ such that $g(p)$ is the
  principal down-set $\doa x$. Then $g_{px}$ is a small map since $\doa x$ is finite.
  Furthermore if
  $h$ is an isotone map such that $h$ is not $\leq f$, then for at least
  one pair $(p,x) \in X$ we must have $x \in h(p)$ and so
  $h \geq g_{px}$. Then the complement
  \[ U(\bfnu,f)^c = \cup_{(p,x) \in X} U(g_{px},\bfen). \]
\end{proof}

\subsection{Criteria for down-sets being regular open or for
  being clopen}
Let
\[ f_1 \leq f_2 \leq f_2 \leq \cdots \]
be a weakly increasing set of maps in $\Pro(P,Q)$, write
$\colim f_r$ for their join. If the $f_r$'s are a decreasing sequence
we write $\lim f_r$ for their meet.

Posets in our application will typically be \napo posets, like finite
posets and the natural numbers $\NN$. We assume in this subsection that
our posets are natural.
The following seems to provide the best way to check whether an open down-set
of isotone maps is regular.
\begin{proposition} \label{pro:wellPQ-close}
  Let $P$ and $Q$ be \napo posets and $\cI$ an {\em open} down-set
  in $\Pro(P,Q)$. 
  \begin{itemize}
  \item[a.] $\cI$ is also closed if and only if it contains the colimit of any increasing
    sequence of isotone maps in $\cI$.
  \item[b.] $\cI$ is regular if and only if it contains any {\em large} colimit of
    an increasing sequence of isotone maps in $\cI$.
  \end{itemize}

  There are analogous statements for open up-sets and decreasing
  sequences of maps in $\cF$. 
\end{proposition}

\begin{proof}
  a1. Suppose $\cI$ is closed. Let $\{f_i\}$ be a weakly increasing sequence
  of isotone maps in $\cI$, and $f = \colim f_i$. If $f$ is {\em not} in $\cI$,
  there  is an open subset $U(b,c)$ in the complement $\cI^c$ containing $f$.

  Since $b(p)$ is finite and $\lim f_i(p) = f(p)$ which contains $b(p)$, there
  is a number $N(p)$ such that $f_i(p)$ contains $b(p)$ for $i \geq N(p)$. 

 Let $T$ be the projection of
 $\Lambda b$ on $P$. Since $b$ is small, $\Lambda b$ is finite and so is $T$.
 Let  $N$ be the maximum of $N(p)$ for $p \in T$. We show below that $f_i \geq
 b$ for $i \geq N$. But then since $f_i \in \cI$ and $\cI$ is a down-set,
 this gives $b \in \cI$, a contradiction since $U(b,c)$ and $\cI$ are disjoint.
 Thus $f$ must be in $\cI$.

 \medskip To show that $f_i \geq b$, if the contrary were true, let
 $q \in b(p)\backslash f_i(p)$ for some $p$ and let $p^\prime \leq p$
 be minimal in $P$ with $q \in b(p^\prime)$. Then:
 \begin{itemize}
\item[i.] $(q,p^\prime) \in \Lambda b$ and so
   $f_i(p^\prime) \supseteq b(p^\prime)$ since $i \geq N \geq N(p^\prime)$. 
 \item[ii.] Since $p^\prime \leq p$ then $f_i(p^\prime) \sus f_i(p)$ and so
   $q \not \in f_i(p^\prime)$. 
 \end{itemize}
 But these give a  contradiction, and so we must have $f_i \geq b$.

 \medskip
\noindent a2. Suppose any limit of an increasing sequence in $\cI$, is in $\cI$.
 Suppose $\cI$ is {\it not} closed. Then there is $f \in \cI^c$
 such that any open subset of $f$ intersects $\cI$. Now there exists the
 following:
 \begin{itemize}
   \item[i)] A sequence of finite posets in $P$
  \[ I^1 \sus I^2 \sus \cdots \]
  such that $\cup_i I^i = P$.
\item[ii)] For each $p \in P$ we can find a sequence
  of finite posets of $Q$
  \[ J^1_p \sus J^2_p \sus \cdots  \]
  such that $\cup_j J_p^j = f(p)$.
\end{itemize}
Givne $k$, define the isotone maps $f_k$ by:
\[  f_k(p) = \cup_{p^\prime \leq p} J^k_{p^\prime}  \text{ when } p \in I^k, \quad
  f_k(p) =  \cup_{p^\prime \in I^k} f_k(p^\prime) \text{ when } p \not \in I^k. \]
Then
\begin{itemize}
\item $f_k$ is isotone and small,
\item $i \leq j$ implies $f_i \leq f_j$,
\item The open subset $U(f_i,\bfen)$ contains $f$ and so intersects
  $\cI$. Then $f_i \in \cI$ since $\cI$ is a down-set,
\item $\colim f_i = f$ and so $f \in \cI$.
\end{itemize}
But the latter contradicts $f \not \in \cI$.

\medskip
b1. Let  $\cI$ be regular and $\{f_i\}$ an increasing sequence in $\cI$, whose
colimit $f$ is large.
By part a the colimit of $f_i$ is in the closure $\overline{\cI}$. Since $\cI$
is regular $\overline{\cI} \backslash \cI$ does not contain large
maps by Lemma \ref{lem:top-gap}. Hence $f \in \cI$.

b2. Suppose $\cI$ contains large limits of increasing sequences in $\cI$.
Let $f$ be a large map in the closure $\overline{\cI}$. By the construction
in a2 there is a sequence $\{ f_i \}$ of small maps such that $f = \colim f_i$.
Each open subset $U(f_i,f)$ intersects $\cI$ and so $f_i$ is in $\cI$.
But then $f$ is in $\cI$. Then by Lemma \ref{lem:top-gap}, $\cI$ is regular.
\end{proof}

\begin{example} \label{eks:monwell-fr}
  Let the open down-set $\cI$ of $\Pro(\NN, \NN)$ be generated by
  the large maps $f_1, f_2, f_3, \cdots$ 
  given by
  \[ f_r(n) = \begin{cases} r & n \leq r \\
      \infty & n > r \end{cases}
    \]
    Then $\colim_r f_r$ is the maximal isotone map $\bfen$, which is not
    in $\cI$. Hence $\cI$ is not regular. Also, its closure is all of
    $\Pro(\NN, \NN)$.
  \end{example}

  \begin{example} \label{eks:monwell-gr}
    Let the open up-set $\cF$ in $\Pro(\NN, \NN)$ be generated by
    the small maps $g_1, g_2, \cdots$ given by 
    \[ g_r(n) = \begin{cases} 1 & n < r \\
        2 & n \geq r
      \end{cases}. \]
    The limit $\lim_r g_r$ is the minimal isotone map $\bfnu$, which is not
    in $\cF$. Hence $\cF$ is not regular. Also, its closure is all of
    $\Pro(\NN, \NN)$.
  \end{example}

  \begin{example} \label{eks:monwell-f}
    Let the profunctor $h : \NN \pil \NN$ be given by
    $h(i) = i$. Let $\cI$ be the down-set generated by the large maps $h^r$ for
    $r \geq 1$, and
    $\cF$ be the up-set generated by the small maps $h_r$ for $r \geq 0$, where
    \[ h^r(i) = \begin{cases} i & i \leq r \\
        r & i = r+1 \\
        \infty & i \geq r+2
        \end{cases}, \quad
h_r(i) = \begin{cases} i & i \leq r \\
        r+2 & i \geq r+1
        \end{cases}. \]
Then $[\cI, \cF]$ is a Dedekind cut and the gap consists of the function $h$.
      \end{example}

      In the argument that follows now we use the following construction:
      Let $I \sus J \sus P$ be down-sets and $f : I \pil \di{Q}$
      an isotone map. Then there exists a unique {\it minimal} extension
      $f^J : J \pil \di{Q}$ such that the restriction$(f^J)_{|I} = f$.
      For each $p \in J$ we let
      \[ f^J(p) = \vee_{p^\prime \in I, p^\prime \leq p} f(p^\prime). \]
      Note that if $f$ is small and $I$ is finite, the extension
      $f^J$ is also small.
      
      The following characterizes precisely when the gap is empty.
\begin{theorem} \label{thm:natural:fg} Let $P$ and $Q$ be natural posets.
  A down-set $\cI$ of $\Pro(P,Q)$ is clopen (closed and open) if and
  only if
  $\cI$ is a finite union of basis open subsets $U(\bfnu,f)$.
  Alternatively formulated, an open down-set $\cI$ is clopen
  if and only if it is finitely generated.
\end{theorem}

\begin{proof}
  The open subsets $U(\bfnu,f)$ are closed, by Proposition
  \ref{pro:wellPQ-close}. Hence any finite
  union of them is clopen.

\medskip
Suppose now $\cI $ is clopen. We will show it has only a finite number
of minimal generators. Assume the opposite, that there are infinitely
many generators of $\cI$, none of which are comparable.

\medskip
\noindent 1. Take a filtration of finite posets
\[ I_1 \sus I_2 \sus \cdots \sus I_r \sus \cdots \]
such that the union of the $I_r$'s is $P$. Let $T_1$ be the set of isotone maps
$F_1 : I_1 \pil \di{Q}$ such that:
\begin{itemize}
  \item[] For every $f : I_1 \pil \di{Q}$ with
    $f \leq F_1$ and the $f(p)$ {\em finite} down-sets in $Q$
    for $p \in I_1$, there are {\it infinitely many} of the 
    minimal generators $g$ of $\cI$ such that the restriction
$g_{|I_1} \geq f$. Clearly $T_1$ is non-empty: we may take $F_1$ to be the
constant map with value the empty set.
\end{itemize}

\medskip
\noindent 2. Let us show that we can use Zorn's lemma on $T_1$. If
\[ F_1^1 \sus F_1^2 \sus \cdots \sus F_1^r \sus \cdots \]
is a chain in $T_1$, let $F_1 = \colim F_1^r$. We claim $F_1$ is
in $T_1$.  Let $f \leq F_1$ with $f$ having finite values. For each $p$,
for some integer
$N(p)$ we have, since $f(p)$ is finite, that $f(p) \sus F^{N(p)}_1(p)$.
Let $N$ be the maximum of
the $N(p)$ for $p \in I_1$. Then $f \leq F_1^N$ and so there are
infinitely many minimal generators $g$ of $\cI$ with $g_{|I_1} \geq f$.
This shows that Zorn's lemma applies.
Furthermore, there is an increasing sequence $\{ f\}$
of such functions whose limit is $F_1$. Extending these $f$'s to $P$
we get an increasing
sequence $\{ f^P \}$ whose limit is the extension $F_1^P$.
Each $f^P$ is in $\cI$ since
there are minimal generators $g$ of $\cI$ with $g_{|I_1} \geq f$.
Since $\cI$ is closed, by Proposition \ref{pro:wellPQ-close}
we have $F_1^P \in \cI$, a fact we use in 5. below. 

\medskip
\noindent 3. Let $F_1$ be maximal in $T_1$. We now show the following refined
statement:
\begin{claim} For each finite $f \leq F_1$ there are infinitely many
  minimal generators $g$ such that $f \leq g_{|I_1} \leq F_1$.
\end{claim}

\begin{proof}
  Let $p_1, p_2, \ldots, p_m$ be the elements of $I_1$ in an order such that
  $p_i < p_j$ implies $i < j$. (This is usually called a linear extension of
  $I_1$.)
  Suppose we have verified that for each $f \leq F_1$ there are infinitely
  many $g$'s such
$f \leq g_{|I_1}$ and with $f(p_i) \leq g(p_i) \leq F_1(p_i)$ for $i < s$.
Let $X_s$ be the (finite) set of minimal elements of the complement
$F_1(p_s)^c$. If not $g(p_s) \sus F_1(p_s)$ then $g(p_s)$ contains some
$x \in X_s$. Define the isotone $F_{1x} : I_1 \pil \di{Q}$ by
\[ F_{1x}(p)  = \begin{cases} F_1(p), & \text{ not } p \geq p_s \\
    F_1(p) \cup \{x \} & p \geq p_s \end{cases}  \]
Since $F_1$ is maximal, there is some $f_x \leq F_{1x}$ where $f_x$ has
finite values, such that there are not infinitely many $g$'s
with $f_x \leq g_{|I_1}$. Let $f_x^\prime$ be the meet $f_x \wedge F_1$. 
Then let $f^\prime$ be the join $f \vee (\vee_{x \in X_s} f_x^\prime)$,
an isotone with finite values. Note that $f^\prime \leq F_1$.
Then there
are infinitely many $g$'s such that: i) $f^\prime \leq g$
and ii) $f^\prime(p_i) \leq g(p_i) \leq F_1(p_i)$ for $i < s$. But only
a finite number of $g$'s  are larger than $f_x$ for each $x$ in the finite
set $X_s$. So taking
away a finite number of the $g$'s, none of the infinitely many
rest will have $g(p_s)$
containing any $x \in X_s$. Hence for these $g(p_s) \sus F_1(p_s)$.
\end{proof}

\medskip
\noindent 4. Suppose we now have constructed $F_j : I_j \pil \di{Q}$ for
$j = 1, \cdots, r$ such that $F_{j|I_i} = F_i$ for each $i \leq j$
and the $F_j$ are such that for each finite $f : I_j \pil \di{Q}$
with $f \leq F_j$
there are infinitely many generators $g$ such that $f \leq g_{|I_j} \leq F_j$,
and $F_j$ is maximal such. We will extend $F_r$ to a function
$F_{r+1} : I_{r+1} \pil \di{Q}$.
Let $T_{r+1}$ be the set of isotone $F : I_{r+1} \pil \di{Q}$ such that
$F_{|I_r} = F_r$ and 
for each finite $f \leq F$ there are infinitely many $g$'s with
$g_{|I_r} \geq f$ and $F_r \geq g_{|I_r} \geq f$. 
The set $T_{r+1}$ is non-empty since it is easy to see that the
extension $F_r^{I_{r+1}}$
fulfills the condition of being in $T_{r+1}$ since $F_r$ does so for
$T_r$. As above $T_{r+1}$ fulfills the condition for Zorn's lemma and
so $T_{r+1}$ has a maximal element $F_{r+1}$. As above we may show that
for each finite $f \leq F_{r+1}$ 
there are infinitely many minimal generators $g$ of $\cI$ such that
$F_{r+1} \geq g_{|I_{r+1}} \geq f$

\medskip
\noindent 5.  Consider now the sequence of extensions
\[ F_1^P \leq F_2^P \leq \cdots \leq F_r^P \leq \cdots . \]
Let $F = \colim F_r^P$. Since each $F_r^P \in \cI$ (by the same reason
as at the end of part 2) and $\cI$ is closed
we get $F \in \cI$. Since $\cI$ is open there exists then 
a large $G \geq F$ in $\cI$. Let $I$ be the finite profile down-set of $G$.
Due to $I$ being finite we must have $I_r \supseteq I$ for sufficiently
large $r$. But there are infinitely many minimal generators $g$ such
that $g_{|I_r} \leq F_{r}$ which again is $\leq G_{|I_r}$.
But since $G(p) = \infty$ for $p$ in the
complement $I_r^c$, we would have $G \geq g$. This is a contradiction since
the $g$'s are (infinitely many) minimal generators of $\cI$.

The upshot is that there cannot be infinitely many generators of
$\cI$.
\end{proof}

\begin{problem}
  What subsets of $\Pro^u(\NN,\NN)$ can be gaps?
\end{problem}
  
\section{Natural posets and finite type cuts}
\label{sec:monwell}
Here we give the variant of Theorem \ref{thm:AD-duality} such
that the $\Lambda$- and $\Gamma$-ideals are finitely generated ideals in
a polynomial ring. This works when $P$ and $Q$ are natural posets,
which we assume in this section.

\medskip
By Lemma \ref{lem:AD-fin}
the map $\Lambda$ of \eqref{eq:AD-Lambda} restricts to a map:
\begin{equation*}
  \Lambda : U\Pro_S(P,Q)  \pil  (UQ \times UP^\op)^\dual_\fin, \quad
  f \mapsto  (\Lambda f,-).
\end{equation*}
Similarly the  map $\Gamma$ of \eqref{eq:AD-Gamma} restricts to a map:
\begin{equation*}
\Gamma : U\Pro^L(P, Q)  \pil  (UQ \times UP^\op)^{\widehat{\fin}}, \quad
  f  \mapsto  (-,\Gamma f).
\end{equation*}
For a cut $(\cI, \cF)$ for $\Pro(P,Q)^\dual$ let $\cF_S$ denote the
small maps in $\cF$, and $\cI^L$ the large maps in $\cI$.

By \eqref{eq:PQ-gu} the map $\Lambda$ induces an isotone map:
\begin{equation*}
  \Lambda_U^\ii : \Pro_S(P, Q)^\dual \pil
  (UQ \times UP^\op)^{\doublehat{}}_{\fin},
  \quad (-,\cF_S)  \mapsto  (-, \Lambda(\cF_S)^\uparrow).
\end{equation*}
The map $\Gamma $ induces an isotone map:
\begin{equation*}
  \Gamma_U^! : \Pro^L(P, Q)^{\dual}  \pil
  (UQ \times UP^\op)^{\doublehat{\fin}},
  \quad (\cI^{L},-)  \mapsto  (\Gamma(\cI^{L})^\downarrow,-).
\end{equation*}
We have the following variation of Theorem \ref{thm:AD-duality}.

\begin{theorem} \label{thm:monwell-duality}
  Let $P$ and $Q$ be natural posets and
  $[\cI,\cF]$ a Dedekind cut for $\Pro(P, Q)$.
  Then $(\Gamma(\cI^L)^{\downarrow},\Lambda(\cF_S)^{\uparrow})$ is a finite type cut for the Boolean lattice
  $(UQ \times UP^\op)^\dual$.
\end{theorem}

\begin{proof}
  We must show 1 and 2 in Definition \ref{def:sett-AD}.
  We do 1 as 2 is analogous.
  
  \noindent {\bf Part I.}
  We show the implication $\Rightarrow$ of 1 in Definition
  \ref{def:sett-AD}.
 This amounts to show
  for any $(I,F) \in  \cI^L $ and $(J,G) \in \cF_S$ that 
  \[ \Gamma F \cap \Lambda J  \neq \emptyset . \]
The argument for this is the same as in Part I for Theorem \ref{thm:AD-duality}.

\medskip
\noindent {\bf Part II.} We show the implication $\Leftarrow$ of 1 in Definition
\ref{def:sett-AD}.
Let $S \sus UQ \times UP^\op$ be a {\em finite} set such that
$S \cap \Lambda J \neq \emptyset$ for {\it every} $(J,G) \in \cF_S$.
We show that $S \supseteq \Gamma F$ for {\it some} $(I,F) \in \cI^L$.

Let $I_0 \sus P$ and $J_0 \sus Q$ be finite down-sets such that
the projections of $S$ are contained respectively in $I_0$ and in $J_0$.
(We know they are finite due to $P$ and $Q$ being down-finite.)
For a large $f$ in $\Pro(P,Q)$ denote by $\Lambda_0 f$ the subset of
$\Lambda f$ lying above $I_0$. 
Let $T$ consist of the large $f$ such that:
\begin{itemize}
\item $f(p) = \begin{cases} \infty, & p \not \in I_0 \\
      \sus J_0, & p \in I_0
      \end{cases}. $ 
  \item $\Lambda_0 f \cap S = \emptyset$.
\end{itemize}
Clearly the minimum element of those $f$ fulfilling the first point, also
fulfills the second, so $T$ is non-empty.
We use Zorn's lemma to show that $T$ has a maximal element.
So let
\[ f_1 \leq f_2 \leq \cdots \leq f_n \leq \cdots \]
be a chain of maps in $T$ and let $f = \colim f_r$.
We claim that $\Lambda_0 f \cap S = \emptyset$. So suppose
$(q,p) \in \Lambda f$, so $ q \in f(p) $ but $q \not \in f(p^\prime)$ for
every $p^\prime < p$. Then $q \in f_i(p)$ for some $i$ and
$q \not \in f_i(p^\prime)$ since $f_i(p^\prime) \sus f(p^\prime)$. Then
$(q,p) \in \Lambda f_i$ and so is not in $S$. Thus $f$ is an upper bound
in $T$ for the chain. By Zorn's lemma $T$ has a maximal element, which we
denote $f$.

\begin{claim} $f$ is in $\cI$.
\end{claim}

\begin{proof} The map $f$ is not in the gap between $\cI$ and $\cF$ since
  it is large. So if $f$ is not in $\cI$ it must be in $\cF$. Since
  $\cF$ is open,  for $B \supseteq J_0$ large enough
$f_B$ given by
\[ f_B(p) = \begin{cases} f(p), & p \in I_0 \\
    B, & p \not \in I_0
  \end{cases}, \]
is in $\cF$. But $\Lambda_0 f_B \sus \Lambda_0 f$, and
so $\Lambda f_B \cap S = \Lambda_0 f_B \cap S$ is the empty set,
which is not so by the requirement on $S$.
Hence $f$ is not in $\cF$ and not in the gap, so it must be in $\cI$.
\end{proof}

\begin{claim}
  $  \Gamma f \sus S$. This finishes the proof by the second line
  of Part II.
\end{claim}

\begin{proof} Suppose not, so there is
  $(q_0,p_0) \in \Gamma f \backslash S$.
Recall $q_0$ is minimal in the complement $f(p_0)^c$. Let
\[ \tilde{f} (p) = \begin{cases} f(p) \cup \{ q_0 \}, &
    p \geq p_0 \\ f(p), & \text{otherwise}.
  \end{cases} \]
We now claim that:
\[ \Lambda_0 \tilde{f} \sus \Lambda_0 f \cup\{ (q_0,p_0)\}. \]
That $(q,p) \in \Lambda_0 \tilde{f}$ means that $p \in I_0$ and 
  $q \in \tilde{f}(p)$, and
  $q \not \in \tilde{f}(p^\prime)$ for $p^\prime < p$.
  \begin{itemize}
  \item[1.]  If not $p \geq p_0$ then $\tilde{f}(p^\prime)  =
    f(p^\prime)$ for all
  values $p^\prime \leq p$. Therefore $(q,p)$ is in $\Lambda_0 f$.
\item[2.]  Let $p \geq p_0$. Then either i. $q \in f(p)$ or ii. $q = q_0$.
  In case i. $(q,p) \in \Lambda f$. In case ii. where
  $q = q_0$, 
  if $p > p_0$ then $(q_0,p)$ would  not be in $\Lambda \tilde{f}$
  since $q_0 \in \tilde{f}(p)$ and $q_0 \in \tilde{f}(p_0)$.
  We are thus left with $(q,p) = (q_0,p_0)$.
\end{itemize}

The upshot is that $\tilde{f}$ fulfills the requirements to be in $T$.
But this contradicts $f$ being maximal in $T$. Hence no such $(q_0,p_0)$
can exist, and so $\Gamma f \sus S$.
\end{proof}
\end{proof}
\medskip

We may now define the $\Lambda$- and $\Gamma$-ideals for (infinite)
natural posets $P$ and $Q$. These ideals then live in
infinite-dimensional polynomial rings. These ideals
will be square-free monomial ideals.

\begin{definition} Let $[\cI,\cF]$ be a Dedekind cut for
  $\Pro(P,Q)$.
  Let $L_\Lambda(\cF)$ the ideal in $k[x_{Q \times P^\op}]$
  generated by the monomials $\prod_{(q,p) \in I} x_{q,p}$
  for $(I,-) \in \Lambda(\cF_S)$.
  Let $L_\Gamma(\cI)$
    be the ideal in $k[x_{P \times Q^\op}]$ generated by the monomials
    $\prod_{(p,q) \in F^\op}x_{p,q}$ for $(F^\op,-)
    \in \Gamma(\cI^L)^{\op}$.
    \end{definition}
    By the theorem above, these ideals are Alexander dual ideals, meaning
that the monomials in $L_\Lambda(\cF)$ are precisely those monomials
with non-trivial common divisor with every monomial in $L_\Gamma(\cI)$, and
vice versa.

\section{Monomial ideals} \label{sec:mon}

For the case $Q = \NN$ we show that open down-sets in
$\Pro(P,\NN)$ induce monomial ideals in the polynomial ring $k[x_P]$.
In particular when $P = \NN$, these monomial ideals are precisely the
class of strongly stable ideals in the infinite-dimensional
polynomial ring $k[x_\NN]$. The results in Section \ref{sec:top}
then enables defining a duality on a distinguished
class of strongly stable ideals, see \cite{Fl-St}.

\medskip
Given a map of sets $A \mto{\alpha} B$. Let
$\Hom_\fin(A,\NN_0)$ be the maps $A \mto{u} \NN_0$ such that
the support $\{ a \, | \, u(a) > 0 \}$ is a finite set.
Then $\alpha $ induces a map
\[ \Hom_\fin(A,\NN_0) \mto{\tilde{\alpha}} \Hom_\fin(B,\NN_0) \]
sending a map $u$ to a map $v : B \pil \NN_0$ where
$v(b) = \sum_{\alpha(a) = b} u(a)$.

Consider the projection $UP \times UQ^\op \mto{\pi_1} UP$, and inclusion
$\boto \mto{i} \NN_0$ sending $0 \mapsto 0$ and $1 \mapsto 1$.
We get a composite
\[ \Hom_\fin(UP \times UQ^\op, \boto) \mto{\tilde{i}}
\Hom_\fin(UP \times UQ^\op, \NN_0) \mto{\tilde{\pi}_1} \Hom_\fin(UP,\NN_0).\]
Given a map $\tau : UP \times UQ^\op \pil \boto$
corresponding to a finite subset
$T \sus UP \times UQ^\op$ (the arguments for which $\tau$ is non-zero),
then $\tilde{\pi} \circ \tilde{i}(\tau)$
is the map sending $p$ to the cardinality of the set of pairs $(p,q) \in T$
with first coordinate $p$.

Note that the set $\Hom_\fin(UP, \NN_0)$ identifies as the set
of monomials $\Mon(x_P)$ in the variables $x_p$ for $p \in P$. 

\begin{lemma} Let $Q = \NN$. 
  The composite $\lambda = \tilde{\pi_1} \circ \tilde{i} \circ \Lambda$:
  \[ \lambda : U\Pro_S(P,\NN) \mto{\Lambda}
\Hom_\fin(UP \times U\NN^\op, \boto) \mto{\tilde{\pi}_1  \circ \tilde{i}} 
    \Hom_{\fin}(UP,\NN_0) \]
  is an isomorphism.
\end{lemma}

\begin{proof}
  This is Proposition 4.3 in \cite{F-LP}.
\end{proof}

We may also consider the composite
$\gamma = \tilde{\pi}_2 \circ \tilde{i} \circ \Gamma$:
\[ \gamma : U\Pro^L(\NN,{P}) \mto{\Gamma}
\Hom_\fin(U\NN \times UP^\op, \boto) \mto{\tilde{\pi}_2 \circ \tilde{i}} 
  \Hom_{\fin}(UP^\op,\NN_0), \]
which is likewise an isomorphism.
\begin{lemma}
  There is a commutative diagram:
  \[ \xymatrix{ U\Pro_S(P,\NN) \ar[rd]^{\lambda} \ar[dd]^D & \\
      &  \Hom_\fin(UP,\NN_0). \\
      UPro^L(\NN,{P})^\op \ar[ru]^{\gamma^{\op}} &
    }
  \]
  Note that the vertical map sends a cut $(I,F)$
  to a cut $(F^\op, I^\op)$.
\end{lemma}

As in \eqref{eq:PQ-gu} this extends to a commutative diagram:

\[ \xymatrix{  \Pro_S(P,\NN)^\dual \ar[rd]^{\lambda^\ii} \ar[dd]^D & \\
      &  \Hom_\fin(UP,\NN_0)^\dual. \\
      \Pro^L(\NN,{P}{\widehat{)^\op}} \ar[ru]^{\gamma^{! \, \op}} &
    }
  \]
  The up-sets in $\Hom_\fin(UP,\NN_0)$
  identifies as the monomial ideals
  in $k[x_P]$.
  Recall that open up-sets $\cF$ in $\Pro(P,\NN)$ are in bijection with
  up-sets in $\Pro_S(P,\NN)$. So we get an upset $\cF_S$ in the latter.
  By $\lambda^{\ii}$ we get a monomial ideal in $k[x_P]$.

  \medskip
\noindent {\bf Application.} When $P = \NN $ the monomial ideals
  one gets are precisely the strongly stable ideals in $k[x_\NN]$.
  There is a one-one correspondence between open up-sets in
  $\Pro(\NN,\NN)$ and strongly stable ideals in $k[x_\NN]$,
  see \cite{Fl-St}.

  \begin{definition}
  Given a Dedekind cut $[\cI,\cF]$ for $\Pro(\NN, \NN)$.
  We get an up-set $\lambda^\ii(\cF_S)$ in $\Hom_\fin(\NN,\NN_0)$
which corresponds to a monomial ideal $I_\lambda$ of $k[x_\NN]$.
    Similarly we get a up-set $\gamma^{! \op}(\cI^{L, \op})$ in
    $\Hom_\fin(\NN,\NN_0)$ which corresponds to a monomial ideal $I_\gamma$ 
of $k[x_\NN]$. The ideals $I_\lambda$ and $I_\gamma$ are {\em dual
strongly stable ideals}.
\end{definition}

Note that this definition is reasonable. Given $I_\lambda$ we 
can reconstruct $\cF_S$ and thereby $\cF$. Since 
$[\cI,\cF]$ is a Dedekind cut, we can construct $\cI$ from $\cF$ and
thereby determine $I_\gamma$. And of course we may go the opposite
way, from $I_\gamma$ we may determine $I_\lambda$. Note also that not
every strongly stable ideal $I$ has a dual. Only the strongly stable ideals
such that the corresponding open down-set $\cI$ is {\em regular}, have
duals. Versions of this duality for finite-dimensional polynomial rings
were first given in \cite[Section 7]{F-LP} and independently by \cite{ShiYan}.

\medskip
{\it Funding.}
The authors did not receive support from any organization for the submitted
work.

\medskip
{\it Competing interests.}
The authors have no competing interests to declare that are relevant to the content of this article.

\bibliographystyle{amsplain}
\bibliography{biblio}

\providecommand{\bysame}{\leavevmode\hbox to3em{\hrulefill}\thinspace}
\providecommand{\MR}{\relax\ifhmode\unskip\space\fi MR }
\providecommand{\MRhref}[2]{%
  \href{http://www.ams.org/mathscinet-getitem?mr=#1}{#2}
}
\providecommand{\href}[2]{#2}
\begin{thebibliography}{10}

\bibitem{Ba}
John~C Baez, \emph{Isbell duality}, Notices of the American Mathematical
  Society \textbf{70} (2022), no.~1.

\bibitem{Ben}
Jean B{\'e}nabou, \emph{Distributors at work}, Lecture notes written by Thomas
  Streicher \textbf{11} (2000).

\bibitem{Bor}
Francis Borceux, \emph{{Handbook of categorical algebra: volume 1, Basic
  category theory}}, vol.~1, Cambridge University Press, 1994.

\bibitem{Ca}
Aurelio Carboni and Ross Street, \emph{Order ideals in categories}, Pacific
  journal of mathematics \textbf{124} (1986), no.~2, 275--288.

\bibitem{DFN}
Alessio D’Al{\`\i}, Gunnar Fl{\o}ystad, and Amin Nematbakhsh,
  \emph{{Resolutions of co-letterplace ideals and generalizations of Bier
  spheres}}, Transactions of the American Mathematical Society \textbf{371}
  (2019), no.~12, 8733--8753.

\bibitem{EbSe}
Carl Eberhart and John Selden, \emph{On the closure of the bicyclic semigroup},
  Transactions of the American Mathematical Society \textbf{144} (1969),
  115--126.

\bibitem{EHM}
Viviana Ene, J{\"u}rgen Herzog, and Fatemeh Mohammadi, \emph{{Monomial ideals
  and toric rings of Hibi type arising from a finite poset}}, European Journal
  of Combinatorics \textbf{32} (2011), no.~3, 404--421.

\bibitem{F-LP}
Gunnar Fl{\o}ystad, \emph{{Poset ideals of P-partitions and generalized
  letterplace and determinantal ideals}}, Acta Mathematica Vietnamica
  \textbf{44} (2019), no.~1, 213--241.

\bibitem{Fl-St}
\bysame, \emph{{Shift modules, strongly stable ideals, and their dualities}},
  arXiv preprint arXiv:2105.14604 (2021), 1--39.

\bibitem{FGH}
Gunnar Fl{\o}ystad, Bj{\o}rn~M{\o}ller Greve, and J{\"u}rgen Herzog,
  \emph{Letterplace and co-letterplace ideals of posets}, Journal of Pure and
  Applied Algebra \textbf{221} (2017), no.~5, 1218--1241.

\bibitem{FV-Bi}
Gunnar Fl{\o}ystad and Jon~Eivind Vatne, \emph{{(Bi-) Cohen--Macaulay
  simplicial complexes and their associated coherent sheaves}}, Communications
  in algebra \textbf{33} (2005), no.~9, 3121--3136.

\bibitem{ACT}
Brendan Fong and David~I Spivak, \emph{An invitation to applied category
  theory: seven sketches in compositionality}, Cambridge University Press,
  2019.

\bibitem{GuSn}
Sema G{\"u}nt{\"u}rk{\"u}n and Andrew Snowden, \emph{The representation theory
  of the increasing monoid}, arXiv preprint arXiv:1812.10242 (2018).

\bibitem{HeQuSh}
J{\"u}rgen Herzog, Ayesha~Asloob Qureshi, and Akihiro Shikama, \emph{Alexander
  duality for monomial ideals associated with isotone maps between posets},
  Journal of Algebra and Its Applications \textbf{15} (2016), no.~05, 1650089.

\bibitem{JuKaMa}
Martina Juhnke-Kubitzke, Lukas Katth{\"a}n, and Sara~Saeedi Madani,
  \emph{Algebraic properties of ideals of poset homomorphisms}, Journal of
  Algebraic Combinatorics \textbf{44} (2016), no.~3, 757--784.

\bibitem{Ju-Survey}
Martina Juhnke-Kubitzke and Sara~Saeedi Madani, \emph{{Ideals Associated to
  Poset Homomorphisms: A Survey}}, Homological and Computational Methods in
  Commutative Algebra, Springer, 2017, pp.~129--140.

\bibitem{Kr}
Joseph~B Kruskal, \emph{{The theory of well-quasi-ordering: A frequently
  discovered concept}}, Journal of Combinatorial Theory, Series A \textbf{13}
  (1972), no.~3, 297--305.

\bibitem{Mi-St}
Ezra Miller and Bernd Sturmfels, \emph{Combinatorial commutative algebra}, vol.
  227, Springer Science \& Business Media, 2004.

\bibitem{NaRo}
Uwe Nagel and Tim R{\"o}mer, \emph{{Equivariant Hilbert series in
  non-noetherian polynomial rings}}, Journal of Algebra \textbf{486} (2017),
  204--245.

\bibitem{CCD4}
Robert Rosebrugh and Richard~J Wood, \emph{Constructive complete distributivity
  iv}, Applied categorical structures \textbf{2} (1994), no.~2, 119--144.

\bibitem{RW04}
\bysame, \emph{Split structures}, Theory Appl. Categ \textbf{13} (2004),
  no.~12, 172--183.

\bibitem{ShiYan}
Kosuke Shibata and Kohji Yanagawa, \emph{Alexander duality for the alternative
  polarizations of strongly stable ideals}, Communications in Algebra
  \textbf{48} (2020), no.~7, 3011--3030.

\end{thebibliography}

\end{document}